\documentclass[a4paper]{article}
\usepackage{amsthm,amsmath,amssymb, bbm}
\usepackage{color}
\newtheorem{teor}{Theorem}[section]
\newtheorem{defin}[teor]{Definition}
\newtheorem{lemm}[teor]{Lemma}
\newtheorem{osse}[teor]{Remark}
\newtheorem{prop}[teor]{Proposition}
\newtheorem{defi}[teor]{Definition}
\newtheorem{coro}[teor]{Corollary}
\newtheorem{prob}[teor]{Problem}
\newtheorem{assu}[teor]{Assumption}

\newcommand{\bele}{\begin{lemm}\begin{sl}}
\newcommand{\enle}{\end{sl}\end{lemm}}
\newcommand{\bedef}{\begin{defi}\begin{sl}}
\newcommand{\eddef}{\end{sl}\end{defi}}
\newcommand{\bete}{\begin{teor}\begin{sl}}
\newcommand{\ente}{\end{sl}\end{teor}}
\newcommand{\beos}{\begin{osse}\begin{rm}}
\newcommand{\eddos}{\end{rm}\end{osse}}
\newcommand{\beas}{\begin{assu}\begin{rm}}
\newcommand{\eddas}{\end{rm}\end{assu}}
\newcommand{\bepr}{\begin{prop}\begin{sl}}
\newcommand{\empr}{\end{sl}\end{prop}}
\newcommand{\bepro}{\begin{prob}\begin{rm}}
\newcommand{\empro}{\end{rm}\end{prob}}
\newcommand{\bede}{\begin{defin}\begin{sl}}
\newcommand{\edde}{\end{sl}\end{defin}}
\newcommand{\beco}{\begin{coro}\begin{sl}}
\newcommand{\enco}{\end{sl}\end{coro}}

\newcommand{\quext}{\quad\text}
\newcommand{\qquext}{\qquad\text}
\newcommand{\de}{\partial}

\newcommand{\RR}{\mathbb{R}}

\newcommand{\NN}{\mathbb{N}}

\newcommand{\beeq}[1]{\begin{equation}\label{#1}}
\newcommand{\eddeq}{\end{equation}}

\newcommand{\beeqa}[1]{\begin{eqnarray}\label{#1}}
\newcommand{\eddeqa}{\end{eqnarray}}

\newcommand{\beal}[1]{\begin{align}\label{#1}}
\newcommand{\eddal}{\end{align}}

\newcommand{\bespl}[1]{\begin{split}\label{#1}}
\newcommand{\edspl}{\end{split}}

\newcommand{\bega}[1]{\begin{gather}\label{#1}}
\newcommand{\edga}{\end{gather}}

\newcommand{\beeqax}{\begin{eqnarray*}}
\newcommand{\eddeqax}{\end{eqnarray*}}

\def\qed{\ifmmode 
  \else \leavevmode\unskip\penalty9999 \hbox{}\nobreak\hfill
  \fi
  \quad\hbox{\hskip.5em\vrule width.4em height.6em depth.05em\hskip.1em}}
\def\endproofsym{\qed}

\def\endnobox{\def\endproofsym{}\end{proof}\def\endproofsym{\qed}}

\newcommand{\no}{\nonumber}

\newcommand{\beeqao}{\begin{eqnarray}\no}
\newcommand{\bealo}{\begin{align}\no}
\newcommand{\besplo}{\begin{split}\no}
\newcommand{\begao}{\begin{gather}\no}

\newcommand{\duav}[1]{\langle{#1}\rangle}

\newcommand{\io}{\int_\Omega}

\newcommand{\OO}{_{\Omega}}

\newcommand{\bn}{\boldsymbol{n}}

\newcommand{\dn}{\partial_{\bn}}
\newcommand{\fhi}{\varphi}
\newcommand{\Fhi}{\Phi}

\newcommand{\lhs}{left-hand side}
\newcommand{\rhs}{right-hand side}

\DeclareMathOperator{\dive}{div}
\DeclareMathOperator{\deriv}{d}

\DeclareMathOperator{\dist}{dist}
\DeclareMathOperator{\dista}{d}
\DeclareMathOperator{\dom}{dom}

\DeclareMathOperator{\Id}{Id}

\DeclareMathOperator{\loc}{loc}

\newcommand{\LDV}{L^2(0,T;V)}

\let\TeXchi\chi
\def\chi{{\setbox0 \hbox{\mathsurround0pt
$\TeXchi$}\hbox{\raise\dp0 \copy0 }}}

\newcommand{\gammaciapo}{\widehat{\gamma}}

\newcommand{\betaciapo}{\widehat{\beta}}

\newcommand{\calA}{{\mathcal A}}

\newcommand{\calE}{{\mathcal E}}

\newcommand{\calN}{{\mathcal N}}

\newcommand{\calB}{{\mathcal B}}

\newcommand{\hgiu}{\underline{h}}
\newcommand{\hsu}{\overline{h}}
\newcommand{\kgiu}{\underline{k}}
\newcommand{\ksu}{\overline{k}}

\newcommand{\dit}{\deriv\!t}

\newcommand{\ov}{\overline}

\newcommand{\ddt}{\frac{\deriv\!{}}{\dit}}

\newenvironment{giuliorev}{\color{red}}{\color{black}}
\newcommand{\III}{\begin{giuliorev}}
\newcommand{\EEE}{\end{giuliorev}}

\numberwithin{equation}{section}
\begin{document}

\title{Global attractor for a Cahn-Hilliard-chemotaxis model with logistic degradation}

\author{Giulio Schimperna\\
Dipartimento di Matematica, Universit\`a di Pavia,\\
and Istituto di Matematica Applicata\\
e Tecnologie Informatiche ``Enrico Magenes'' (IMATI),\\
Via Ferrata~5, I-27100 Pavia, Italy\\
E-mail: {\tt giulio.schimperna@unipv.it}
\and
Antonio Segatti\\
Dipartimento di Matematica, Universit\`a di Pavia,\\
and Istituto di Matematica Applicata\\
e Tecnologie Informatiche ``Enrico Magenes'' (IMATI), \\
Via Ferrata~5, I-27100 Pavia, Italy\\
E-mail: {\tt antonio.segatti@unipv.it}
}


\maketitle
\begin{abstract}
 We consider a mathematical model coupling the Cahn-Hilliard 
 system for phase separation with an additional equation 
 describing the diffusion process of a chemical quantity 
 whose concentration influences the physical process.
 The main application of the model refers to tumor progression,
 where the phase variable $\fhi$ denotes the local proportion of 
 active cancer cells and the chemical concentration $\sigma$ may refer to
 a nutrient transported by the blood flow or to a drug administered to the patient.
 The resulting system is characterized by cross-diffusion effects
 similar to those appearing in the Keller-Segel model for chemotaxis; 
 in particular, the nutrient tends to be attracted towards the regions where 
 more active tumor cells are present (and consume it in a quickier way).
 Complementing various recent results on related models, we investigate here
 the long-time behavior of solutions under the perspective of 
 infinite-dimensional dynamical systems. To this aim, we first 
 identify a regularity setting in which the system is well posed and 
 generates a closed semigroup according to the terminology introduced
 by Pata and Zelik. Then, partly based on the approach introduced 
 by Rocca and the first author for the Cahn-Hilliard system
 with singular potential, we prove that the semigroup is strongly  
 dissipative and asymptotically compact so guaranteeing the existence
 of the global attractor in a suitable phase space. Finally, we 
 discuss the sign properties of $\sigma$ and show that, 
 if the initial datum $\sigma_0$ is a.e.~strictly positive in 
 the reference domain $\Omega$ with $\ln \sigma_0\in L^1(\Omega)$,
 then, for every $0<\tau<T<\infty$, there exists $\delta=\delta(\tau,T)>0$ such that  
 $\sigma (\cdot,\cdot)\ge \delta>0$ a.e.~in $\Omega\times(\tau,T)$. 
 It is not clear, however, whether this property is uniform in $T$;
 actually, we cannot exclude that $\delta$ could degenerate to $0$ 
 as $T$ is let go to infinity.
\end{abstract}


\vspace{2mm}

\noindent {\bf AMS (MOS) subject clas\-si\-fi\-ca\-tion:}%
~~35B41, 35D30, 35K35, 35Q92, 92C17.

\bigskip
\bigskip
\smallskip


\section{Introduction}
\label{sec:intro}

Letting $\Omega\subset \RR^d$, $d\in\{2,3\}$, be a smooth bounded domain of 
boundary $\Gamma:=\de\Omega$, we consider the following PDE system:
\begin{align}\label{CH1}
  & \fhi_t - \Delta \mu = 0,\\
 \label{CH2}
  & \mu = - \Delta \fhi + f(\fhi) - \chi \sigma,\\
 \label{nutr}
  & \sigma_t - \Delta\sigma + \chi \dive( \sigma \nabla \fhi )
    = - h(\sigma,\fhi) \sigma^2 + k(\sigma,\fhi) \sigma,
\end{align} 
where the unknown variables $\fhi$ (the normalized order parameter of a phase 
separation process), $\mu$ (the associated chemical potential) and $\sigma$
(the concentration of a chemical substance) are defined
for $x\in \Omega$ and $t \ge 0$. In the above model, relations \eqref{CH1}-\eqref{CH2}
constitute a suitable form of the Cahn-Hilliard system, where $f$ is 
the derivative of a configuration potential $F$ (described below),
whereas the second-order relation \eqref{nutr} is characterized by
a cross-diffusion term with the same structure as in the Keller-Segel
model for chemotaxis. As first observed in \cite{RSS}, this choice 
prescribes the chemical substance represented by $\sigma$ to migrate towards the 
regions where one of the two phases is prevailing, in a way that is proportional
to its local concentration. The magnitude of this effect is determined
by the {\sl chemotactic coefficient}\/ $\chi>0$ (and, of course, it
is larger for larger $\chi$). The logistic forcing term on the 
\rhs\ of \eqref{nutr}, where the functions $h,k$ are assumed to 
be bounded, strictly positive, and Lipschitz continuous, from the 
one side has a sign-preserving effect guaranteeing (at least)
the nonnegativity of $\sigma$, and, from the other side, provides 
a stabilizing effect as the quadratic term acts a source of 
a-priori estimates preventing the blow-up of $\sigma$. We 
note, however, that this effect is in competition with 
the quadratic growth of the cross-diffusion term, so that the 
actual behavior of solutions somehow depends on a comparison
between the magnitude of the logistic source and that of the chemotactic
coefficient $\chi$, see Remark~\ref{oss:compa} below for more details.

The system is complemented by no-flux (i.e., homogeneous Neumann)
boundary conditions for all unknowns. Concerning $\mu$ and $\sigma$, 
such conditions ensure, respectively, the conservation of mass 
and the absence of inflow (or outflow) of the chemical
through $\Gamma$. Concerning $\fhi$, the no-flux condition 
in \eqref{CH2} corresponds to an orthogonality property of the 
(diffuse) interface with respect to $\Gamma$, a feature
which is standardly assumed in Cahn-Hilliard-based models.

System \eqref{CH1}-\eqref{nutr} is mainly motivated by its applications
to cancer growth processes. In such a setting, $\fhi$ represents
the local proportion of active tumor cells, whereas $\sigma$ denotes the 
concentration or a nutrient (e.g., blood, or some protein, or possibly
also some drug administered to the patient) affecting the evolution 
of the tumor mass. Actually, diffuse interface models for cancer 
evolution are becoming increasingly popular in the recent mathematical
literature. In particular, many recent papers (see, e.g., 
\cite{CL, CLLW, FLRS, GL1, GL2, WLFC} 
and the references therein) have been devoted to the situation
where \eqref{nutr} is replaced by a different relation, namely 
\begin{equation}\label{nutrlin}
  \sigma_t - \Delta\sigma + \chi \Delta \fhi 
    = b(\fhi,\sigma).
\end{equation}
Compared to \eqref{nutr}, \eqref{nutrlin} is mathematically simpler 
as it does not account for quadratic cross-diffusion effects. 
On the one hand, this prevents the 
blow-up of solutions even when the forcing term $b(\fhi,\sigma)$ 
is not of logistic type (but, for instance, has linear growth),
but, on the other hand, relation \eqref{nutrlin} has no sign-preserving 
properties. Consequently, one cannot avoid that $\sigma$ may attain 
negative values (which are somehow ``nonphysical'' as $\sigma$ 
represents a concentration) even if the initial datum
$\sigma_0\ge 0$ a.e.~in~$\Omega$. This fact, as first observed in 
\cite{RSS}, seems to be the main physical motivation for rather 
preferring relation \eqref{nutr}.

As mentioned above, the function $f$ in \eqref{CH2} denotes the derivative,
or, more precisely, the subdifferential, of a configuration potential $F$
of {\sl singular}\/ type, meaning that $F$ is assumed to take 
finite values only in the interval $[-1,1]$ (or, possibly,
only in $(-1,1)$), and to be identically $+\infty$ outside
that interval, with the most relevant choice being represented 
by the so-called Flory-Huggins logarithmic potential defined by
\begin{equation}\label{Flog}
  F(r) = (1+r)\log(1+r)+(1-r)\log(1-r)-\frac\lambda2 r^2, \quad r\in[-1,1],
  \quad \lambda \ge 0.
\end{equation}
Actually, as is customary for phase-separation models, the order
parameter is normalized in such a way that the pure configurations are represented
by $\fhi=\pm 1$; correspondingly, the minima of $F$ are attained in 
proximity of these states, and they are deeper for larger $\lambda$.
Correspondingly, the ``unphysical'' states lying outside $[-1,1]$ 
are penalized by letting $F(\fhi)\equiv +\infty$ for $|\fhi|>1$
(see Assumption~(A2) below for details). We point out that, for 
this class of models, the choice of a ``singular'' potential 
like \eqref{Flog} is furtherly motivated by the fact that it guarantees
the a-priori uniform boundedness of $\fhi$ and, in turn, the 
coercivity of the energy functional (cf.~\eqref{energy} below),
which would be lost in the case when, e.g., $F$ is a ``smooth''
double well polynomial potential (e.g.\ given by $F(r)=(r^2-1)^2$,
which is a  popular choice in the Cahn-Hilliard literature).
To be precise, here we will consider a general class of singular
potentials (including \eqref{Flog} as a particular case),
so to apply to the present model (and in particular to its asymptotic 
analysis for large times) the abstract approach devised in 
\cite{RS} which does not rely on the specific expression 
\eqref{Flog}.

As already observed, the first attempt in considering the coupling
between the Cahn-Hilliard system \eqref{CH1}-\eqref{CH2}
and the Keller-Segel-like relation
\eqref{nutr} was carried out in \cite{RSS}, where the terminology
``Cahn-Hilliard-Keller-Segel'' model was proposed for this class 
of PDE systems. Since then, several papers have been 
devoted to analyzing this model or variants of it, see, e.g., 
\cite{AS,CGSS,GSS,GHW,LRS,RSS,Ssens} and the references therein.
In these papers, existence and regularity properties of 
solutions have been addressed in several interesting situations
(two- or three- space dimensions, occurrence or not of logistic
forcing in \eqref{nutr} and of mass source in \eqref{CH1}, 
presence of nonlinear sensitivity, different types of 
potentials, regularization of approximation schemes).
Moreover, various different solution concepts 
have been proposed, including \cite{LRS} ``entropic''
solutions, whose global existence can be proved even in three space
dimensions with no need for stabilizing terms (a situation
that is very different compared to the genuine Keller-Segel
case for which blow-up occurs in that setting).

A closer look at this recent literature suggests in particular that, from the 
mathematical viewpoint, system \eqref{CH1}-\eqref{nutr} presents some 
specific features that make it rather different both from other
Cahn-Hilliard-based models and from ``standard'' models 
of Keller-Segel type. Actually, on the one hand, 
the poor summability of the forcing 
term $\chi\sigma$ in \eqref{CH2} lowers the regularity properties 
of $\fhi$ (in our case, however, the logistic source 
gives a relevant help) compared with what standardly happens for
Cahn-Hilliard; on the other hand, $\fhi$ has a fourth, rather 
than second order, dynamics with respect to space variables, which 
is something unusual in the Keller-Segel setting.

In this paper, we complement previous mathematical results
on the Cahn-Hilliard-Keller-Segel model by analyzing the long-time behavior
of solutions within the approach of infinite-dimensional 
dynamical systems and with the main purpose of proving existence
of the global attractor in a suitable functional setting.
This seems, indeed, to be the first attempt to address the 
long-time behavior of solutions to this model as previous works seem to
consider only the case of finite time domains.
As said, we deal with the case when a logistic forcing term 
in present in \eqref{nutr} and no mass source occurs in \eqref{CH1}:
in such a setting, it is not too difficult to identify a regularity
class in which well-posedness can be shown. Actually,
in order to have uniqueness, we need to prove a contractive estimate,
and this requires rather good properties both for $\fhi$ and $\sigma$
in order to be able to deal with the quadratic cross-diffusion
term in \eqref{nutr}. In particular, for what concerns 
$\fhi$, the regularity we need (cf.\ \eqref{defi:Fhim} below)
corresponds to the so-called ``energy'' class (according to the 
Cahn-Hilliard terminology), whereas for $\sigma$ the minimal 
condition we are able to consider is a regularity setting of 
``fractional'' type (cf.\ \eqref{defi:Sigma} below). We also
notice that, in view of the occurrence of the singular potential,
the phase-space is endowed with a nonstandard metric (cf.\ \eqref{dist:Fhim}), 
which is mutuated from the approach introduced in \cite{RS} for
the Cahn-Hilliard system with singular
potential. Then, assuming some compatibility condition on coefficients
(cf.\ \eqref{keycoeff} and Remark~\ref{oss:compa} below), we can prove 
that the system is well-posed and generates a closed (according to
the terminology in \cite{PZ}) and 
asymptotically compact semigroup, which in turn
guarantees the existence of the global attractor. 
As in \cite{RSS}, the proof of 
asymptotic compactness is based on the so-called ``second energy
estimate'' for the Cahn-Hilliard system, a procedure providing 
(supposedly) optimal regularity properties for $\fhi$
compatibly with the presence of a singular potential.

After proving existence of the global attractor, we investigate
the sign properties of $\sigma$. Indeed, as relation 
\eqref{nutr} is evidently sign-preserving, the choice \eqref{defi:Sigma}
of the phase space for $\sigma$ needs to incorporate this information.
On the other hand, it is also not difficult to see that,
if one assumes stronger positivity conditions on the initial
datum $\sigma_0$ (for instance, the assumption $\sigma_0>0$
almost everywhere with $\ln\sigma_0\in L^1(\Omega)$ used in \cite{LRS}
for the construction of the ``entropic solutions''),
then this piece of information is maintained
at least on finite time intervals. Hence, it may be worth investigating 
whether, by applying parabolic regularization properties,
one could prove that $\sigma$ is separated from $0$ in
the uniform norm for every value of the time variable $t>0$.
Actually, we can show that the ``entropic'' reformulation
of \eqref{nutr} proposed in \cite{LRS} is suitable for applying
a Moser iteration argument adapted from \cite{SSZ}, which implies
that $\sigma(x,t)\ge \delta(\tau,T)>0$ for a.e.\
$(x,t)\in\Omega\times (\tau,T)$, where $0<\tau<T$.
However, as the estimates 
used in the proof apparently lack dissipativity, we cannot
exclude that the value of $\delta$ might degenerate to $0$
if one lets $T\nearrow \infty$ (see Remark~\ref{oss:nodiss}
below for additional considerations on this point).

\smallskip

\noindent%
The plan of the paper is as follows: in the next section we introduce our 
assumptions in a rigorous framework and state all our mathematical results,
whose proofs are detailed in the remainder of the paper. More precisely,
Section~\ref{sec:apriori} contains the a-priori estimates necessary
for proving existence of weak solutions and dissipativity of the 
associated semigroup; there, the procedure to pass to the limit
in a (hypothetical) approximation scheme is also sketched.
The proof of the well-posedness of the system 
is completed in Section~\ref{sec:uni}, which is devoted to uniqueness,
whereas the subsequent Section~\ref{sec:asy} deals with the parabolic
regularization estimates entailing asymptotic compactness of the solution
semigroup (and consequently existence of the global attractor). Finally,
in Section~\ref{subsec:sig} a reinforced minimum principle property for 
equation~\eqref{nutr} is shown.


\section{Assumptions and main results}
\label{sec:main}

We start with introducing a set of notation which will be useful in
order to rigorously formulate our mathematical results. 
Letting $\Omega$ be a smooth bounded domain of $\RR^d$, 
$d\in\{2,3\}$, of boundary $\Gamma$, we
set $H := L^2(\Omega)$ and $V := H^1(\Omega)$.
For the sake of simplicity, in all proofs we will 
directly consider the case when $d=3$ (particularly, when 
dealing with embedding exponents or related inequalities).
Of course, for $d=2$ some of the results might be improved.

We will often write $H$ in place of $H^d$ 
(with similar notation for other spaces), in case vector-valued 
functions are considered. We denote by $(\cdot,\cdot)$ the 
standard scalar product of~$H$ and by 
$\| \cdot \|$ the associated Hilbert norm.  Moreover, we equip $V$ 
with the standard norm $\|\cdot\|_V^2 = \|\cdot\|^2 + \|\nabla \cdot\|^2$.
Identifying $H$ with its dual space $H'$ by means of the 
scalar product introduced
above, we obtain the chain of continuous 
and dense embeddings $V\subset H \subset V'$.
We indicate by $\duav{\cdot,\cdot}$
the duality pairing between $V'$ and $V$, or, more generally,
between $X'$ and $X$ where $X$ is a generic Banach space continuously
and densely embedded into $H$. Letting $\bn$ denote the outer
unit normal vector to $\Omega$, we also set 
\begin{equation}\label{defi:W}
  W:=\{ v \in H^2(\Omega):~\dn v = 0~\text{on }\, \partial \Omega\},
\end{equation}
which is a closed subspace of $H^2(\Omega)$ and in particular inherits
its norm.

Next, for any functional $\xi\in V'$, we note as 
\begin{equation}\label{media}
  \xi\OO:= \frac1{|\Omega|} \duav{\xi,1}
\end{equation}
the spatial average of $\xi$, where the duality can be replaced by
an integral in case, say, $\xi\in H$. We note as $H_0$, $V_0$ and 
$V_0'$ the closed subspaces of $H$, $V$, $V'$, respectively,
consisting of the function(al)s having zero spatial average. 
Moreover, we observe that the Neumann Laplacian 
$B:=(-\Delta)$, interpreted as a bounded linear operator
\begin{equation}\label{defi:B}
  B: V \to V', \qquad
   \duav{Bv,w}:= \io \nabla v\cdot \nabla w,
\end{equation}
for $v,w\in V$, takes values in $V_0'$ and is invertible as it is
restricted to the functions $v\in V_0$. We shall denote 
as $\calN:V_0' \to V_0$ the inverse of $B$ acting over the 
functionals with zero spatial mean. Then, we observe that 
the norm 
\begin{equation}\label{nor:*}
  \| \xi \|_*^2 := \duav{\xi,\calN \xi}, \qquad
   \xi \in V_0',
\end{equation}
is a norm on $V_0'$ which is equivalent to the standard 
(dual) norm inherited from $V'$. We will use the
above norm on occurrence. In particular, we may notice that,
for $\xi\in V_0'$,
\begin{equation}\label{nor:*2}
   \io | \nabla \calN \xi |^2 
    = \duav{B \calN \xi, \calN \xi}
    = \duav{\xi, \calN \xi} = \| \xi \|_*^2.
\end{equation}

\smallskip
\noindent%
{\bf Assumption (A1).}~~%
We assume the functions $h$ and $k$ on the \rhs\ of \eqref{nutr} to satisfy 
the following structure hypotheses:
\begin{align}\label{hk:reg}
  & h,k\in C^{1}([0,+\infty)\times[-1,1]), \qquad 
   \nabla h, \nabla k\in L^\infty([0,+\infty)\times[-1,1]),\\   
 \label{hk:sign}
  & 0 < \hgiu \le h(\sigma,\fhi) \le \hsu, \quad 
    0 < \kgiu \le k(\sigma,\fhi) \le \ksu, \quext{for every }\,(\sigma,\fhi)\in [0,+\infty)\times[-1,1],
\end{align} 
where $\hgiu,\hsu,\kgiu,\ksu$ are positive constants. Roughly speaking,
the effect of $h$ is to prevent the blowup
of $\sigma$ while the effect of $k$ is that of avoiding it to degenerate to $0$. 

\smallskip
\noindent%
{\bf Assumption (A2).}~~%
As we would like to consider a somehow abstract setting, 
we will assume general structure hypotheses on the singular potential $F$ 
along the lines of the approach devised in~\cite{RS}. 
To introduce the latter, we start with considering a (possibly
multivalued) maximal monotone graph $\beta\subset\RR\times \RR$
(we refer to the monographs \cite{Ba,Br} for the underlying
convex analysis background), and, to avoid technicalities, we 
assume with no loss of generality the normalization
\begin{equation}\label{dombeta}
  \overline{\dom\beta}=[-1,1], \qquad 0\in\beta(0),
\end{equation}
where we recall that the {\sl domain}\/ $\dom\beta\subset\RR$ is the set 
of those $r\in\RR$ such that $\beta(r)$ is not empty. The above choice
is consistent with the ansatz that $\fhi=\pm1$ corresponds to the pure
states or configurations. 
Then, we observe (again, see \cite{Ba} or \cite{Br})
that there exists a (unique) convex and lower semicontinuous
function $\betaciapo:\RR\to[0,+\infty]$ with $\overline{\dom\betaciapo}=[-1,1]$
and  $0=\betaciapo(0)$, such that $\beta=\de \betaciapo$,
$\de$ denoting here the subdifferential of convex analysis. We recall
that, by definition, the domain $\dom\betaciapo$ is the set where $\betaciapo$ 
takes finite values. Then, we may ``reconstruct'' $F$ by setting
\begin{equation}\label{defi:fF}
  F(r) := \betaciapo(r) - \frac\lambda2 r^2, \qquad 
   f(r) = \beta(r) - \lambda r,
\end{equation}
where $\lambda\ge 0$ is a given constant. In this way, $F$ 
($f$, respectively) is decomposed into its convex part $\betaciapo$
(monotone part $\beta$, respectively) part and a
quadratic (linear, respectively) perturbation.
We may notice that, in the specific case of the Flory-Huggins potential
\eqref{Flog}, it turns out that $\dom\betaciapo=[-1,1]$
as $\beta(r)=\log(1+r)-\log(1-r)$ is summable for $|r|\sim 1$.
\beos\label{oss:multi}
 As said, in agreement with the general theory of maximal monotone operators, 
 the graph $\beta$ may be multivalued; namely, for some $r\in[-1,1]$
 the (convex) set $\beta(r)$ may contain more than one element. 
 For this reason, equation \eqref{CH2} should rather be formulated as 
 \begin{equation}\label{CH2-multi}
   \mu = - \Delta\fhi + \eta - \lambda\fhi - \chi\sigma,
    \quext{where }\,\eta(\cdot,\cdot) \in \beta(\fhi(\cdot,\cdot))~~
    \text{almost everywhere},
 \end{equation}
 meaning that $\eta$ is a suitable (and measurable)
 {\sl section}\/ of $\beta(\fhi)$. In particular,
 all the forthcoming estimates, as well as the regularity conditions 
 (like, e.g., \eqref{rego:ffhi} below), should in fact be written in
 terms of the section $\eta$. However, to reduce technicalities,
 in the sequel we shall generally treat $\beta$ as a single-valued 
 function. Just in the definition of the ``attractor space'' 
 (see \eqref{defi:Psim} below), we preferred to maintain the notation
 of \cite{RS}, from which this approach is inspired, 
 where the multivalued structure of $\beta$ explicitly occurs.
\eddos
\smallskip
\noindent%
{\bf Phase space.}~~%
In order to address the long-time analysis of system \eqref{CH1}-\eqref{nutr}
by means of the theory of infinite-dimensional dynamical systems, 
we need to identify a suitable functional set where the system 
generates a dissipative process. In particular, 
according to the general theory (see the monographs \cite{Ha,Te}),
for initial data lying in such {\sl phase space}, we need to be able
to prove well-posedness of the system, continuity properties of 
trajectories with respect to the natural topology of the space,
and existence of a bounded absorbing set (i.e., dissipativity).

To build the phase space, we consider the two variables separately:
for what concerns $\fhi$, it is natural to assume the regularity 
corresponding to  the finiteness of the physical energy and to the conservation
of the total mass. Namely, we are naturally led to consider the set
\begin{equation}\label{defi:Fhim}
  \Phi_m:=\big\{ \fhi \in V:~\betaciapo(\fhi) \in L^1(\Omega),~\fhi\OO=m\big\},
\end{equation}
where $m\in(-1,1)$ is the fixed (and conserved in time) mean value of the initial 
datum.

Then, by \cite[Lemma 3.8]{RS}, it turns out that $\Phi_m$ is a complete 
metric space with respect to the metric 
\begin{equation}\label{dist:Fhim}
  \dist(\fhi_1,\fhi_2) 
   := \| \fhi_1 - \fhi_2 \|_V
    + \| \betaciapo(\fhi_1) - \betaciapo(\fhi_2) \|_{L^1(\Omega)}.
\end{equation}
\beos\label{rem:casolog}
In the case of the Flory-Huggins potential \eqref{Flog} (as well as in all 
other cases when $\dom\betaciapo = [-1,1]$), it is apparent that 
the condition $\betaciapo(\fhi)\in L^1(\Omega)$
is equivalent to $-1\le \fhi\le 1$ 
a.e.\ in $\Omega$. Hence, in this case one may simply take 
$\Phi_m$ as the set of those functions $\fhi\in V$ such that 
$\fhi\OO=m$ and $-1\le\fhi\le 1$ almost everywhere.
Moreover, as this is a closed subset of $V$, one can also omit 
the second summand in \eqref{dist:Fhim} and use simply the $V$-norm.
\eddos
\noindent%
For what concerns $\sigma$, the main issue is related to uniqueness.
Actually, as noted in \cite{RSS}, due to the quadratic growth
of the cross-diffusion term this is a nontrivial property when
one considers \eqref{nutr}. In our setting, the minimal 
regularity under which we can prove 
uniqueness (for the whole system) is of ``fractional type''
(and is a bit weaker than the assumptions taken 
in \cite[Theorem~2.8]{RSS} as 
here we do not have a mass source in \eqref{CH1}). 

To introduce it, also inspired by the approach 
used in \cite{GP}, we first define the linear unbounded 
abstract operator 
\begin{equation}\label{defi:A}
  A:= (I - \Delta) :H \to H,   
   \quad D(A)= W,
\end{equation}
where the space $W$, i.e.\ the {\sl domain}\/ of $A$, 
is defined in \eqref{defi:W}. Then, noting that 
$A$ is strictly positive, we can consider the fractional powers $A^s$ 
of $A$ for any $s\in\RR$. We then set 
\begin{equation}\label{defi:Sigma}
  \Sigma:=\big\{ \sigma\in D(A^{1/4})=H^{1/2}(\Omega):~
   \sigma\ge 0~\text{a.e.~in}~\Omega\big\}.
\end{equation}
This is, of course a closed subset of $D(A^{1/4})$ (hence, a complete 
metric space). The phase space for our asymptotic analysis 
is then defined as the product space $\Phi_m\times\Sigma$.
\beos\label{rem:posi}
Condition \eqref{defi:Sigma} accounts in particular for the 
nonnegativity of $\sigma$. Actually, the structure of the cross-diffusion
term in \eqref{nutr} guarantees (the simplest proof of this fact 
is that based on Stampacchia's truncation method) that, if 
$\sigma_0(\cdot)\ge 0$ almost everywhere in $\Omega$, then it 
also holds that $\sigma(\cdot,\cdot)\ge 0$
almost everywhere in $\Omega\times(0,+\infty)$. 
This property will be improved in Theorem~\ref{teo:ln} below.
\eddos

\smallskip 
\noindent%
We are now ready to state our first theorem regarding well-posedness of 
the initial-boundary value problem for system \eqref{CH1}-\eqref{nutr}
for initial data lying in the product space $\Phi_m\times\Sigma$, as
well as dissipativity of the associated dynamical process in the
same functional setting.
\bete\label{teo:wellpo}
 Let assumptions~{\rm (A1)-(A2) hold}. 
 Let, moreover, the following compatibility condition hold:
 \begin{equation}\label{keycoeff}
   \chi^2 \le 3\hgiu.
 \end{equation}
 For some $m\in (-1,1)$, let also
 \begin{equation}\label{hp:init}
    (\fhi_0,\sigma_0) \in \Phi_m\times\Sigma.
 \end{equation}
 Then, there exists one and only one triple $(\fhi,\mu,\sigma)$ defined
 over $\Omega\times(0,\infty)$ such that, for every $T>0$, 
 \begin{align}\label{rego:fhi}
   & \fhi\in H^1(0,T;V')\cap C^0([0,T];V) \cap L^4(0,T;W) \cap L^2(0,T;W^{2,6}(\Omega)),\\
  \label{rego:ffhi}
   & \betaciapo(\fhi)\in C^0([0,T];L^1(\Omega)), \qquad 
    \beta(\fhi) \in L^2(0,T;L^{6}(\Omega)),\\
  \label{rego:mu}
   & \mu\in \LDV,\\
  \label{rego:sigma}
   & \sigma \in H^1(0,T;D(A^{-1/4})) \cap C^0([0,T];D(A^{1/4})) \cap L^2(0,T;D(A^{3/4})),
 \end{align}
 with $\sigma\ge 0$ almost everywhere in $\Omega\times(0,T)$. The triple 
 $(\fhi,\mu,\sigma)$ satisfies equation \eqref{CH2} a.e.~in $\Omega\times(0,T)$
 with the boundary condition 
 \begin{equation}\label{neum:fhi}
    \dn\fhi=0, \qquext{a.e.~on }\,\Gamma\times(0,T),
 \end{equation}
 together with the following weak formulations of\/ \eqref{CH1} and\/ \eqref{nutr}:
 \begin{align}\label{CH1:w}
   & \duav{\fhi_t,\xi} + \io \nabla\mu \cdot\nabla\xi = 0,
    \qquext{a.e.~in }\,(0,T),\\
  \label{nutr:w}
   & \duav{\sigma_t,\xi} + \io \nabla\sigma\cdot \nabla\xi
    -\chi \io \sigma\nabla\fhi\cdot\nabla\xi
   = \io \big(- h(\sigma,\fhi) \sigma^2 + k(\sigma,\fhi) \sigma\big)\xi,
    \qquext{a.e.~in }\,(0,T),
 \end{align} 
 where (in both relations) $\xi\in V$ is a generic test function,
 and the initial conditions 
 \begin{equation}\label{init}
   \fhi|_{t=0} = \fhi_0, \qquad  \sigma|_{t=0} = \sigma_0.
 \end{equation}
 Moreover, if $S(t)$, $t\ge 0$, denotes the semigroup mapping associated 
 to the above class of weak solutions, the dynamical process described by
 $S$ is dissipative in the phase space $\Phi_m\times\Sigma$.
 Namely, there exists a number $R_1\ge 0$ independent of the initial data
 such that the set
 \begin{equation}\label{asso}
    \calB_0:=\big\{ (\phi,s) \in \Phi_m\times\Sigma:~
     \dist(\phi,0) + \| s \|_{D(A^{1/4})} \le R_1\big\}
 \end{equation}
 is\/ {\rm absorbing}, i.e., for any bounded set $B\in \Phi_m\times\Sigma$
 of initial data there exists a time $T_1>0$ depending only on the ``radius''
 $R$ of $B$ in the metric\/ \eqref{dist:Fhim}, such that 
 \begin{equation}\label{asso:2}
    S(t) B \subset \calB_0
     \quext{for every }\,t\ge T_1.
 \end{equation}
\ente 
\noindent%
\beos\label{oss:compa}
 It is worth making some comment on condition \eqref{keycoeff}, which prescribes that 
 the chemotactic effect is ``small'' compared to the stabilizing effect of the 
 logistic term. Actually, existence of weak solutions without assuming 
 this type of condition is still an open problem, at least in our case
 of a {\sl quadratic}\/ logistic source, whereas, as noted in \cite{RSS},
 if the stabilizing part of the logistic term is superquadratic,
 then \eqref{keycoeff} can be avoided. Compared to \cite[Theorem~2.2]{RSS}, 
 here we have been able to improve (i.e.~weaken) a bit condition
 \eqref{keycoeff}, but we still need a form of it in order to close our 
 a priori estimates.
\eddos
\noindent%
As we have introduced the one-parameter family of mappings 
$t\mapsto S(t)$, it is worth noting that $S(\cdot)$
satisfies the semigroup properties in a proper way. In particular,
it is clear that $S(0) = \Id$ (the identity mapping) and 
that there holds the ``concatenation property''
$S(t_1)\circ S(t_2) = S(t_1+t_2)$ for every $t_1,t_2\ge 0$.
Concerning the continuity properties of trajectories, one can first 
notice that, thanks to the continuity conditions in \eqref{rego:fhi}, 
\eqref{rego:ffhi} and \eqref{rego:sigma}, it turns out that 
the mapping $t\mapsto S(t) (\fhi_0,\sigma_0)$ is continuous 
from $[0,\infty)$ to $\Phi_m\times \Sigma$ (endowed with the 
distance \eqref{dist:Fhim}). Regarding continuity with respect 
to the initial data, this can be proved in a rather straighforward
way; it is however even quickier to observe that the contractive
estimate in the uniqueness proof (cf.\ \eqref{co:56d} below) immediately
implies that $S(\cdot)$ is a {\sl closed}\/ semigroup according 
to the terminology introduced by Pata and Zelik in \cite{PZ},
which is sufficient in order to apply all the machinery from the theory
of infinite-dimensional dynamical systems.
It is worth recalling here some 
basic definition from that theory: first of all, if $X$ is the phase space of 
a dynamical process $S(\cdot)$, a set $K \subset X$ is said to be uniformly 
attracting whenever, for any bounded set $B \subset X$, one has 
\begin{equation}\label{def:attracting}
  \lim_{t\nearrow\infty} \dist(S(t)B, K) = 0,
\end{equation}
where $\dist$ denotes the {\sl unilateral Hausdorff distance}\/ 
of the set $S(t)B$ from $K$, with respect to the metric of $X$, 
namely 
\begin{equation}\label{def:Hau}
  \dist(S(t)B, K) := \sup_{y\in S(t)B}
   \inf_{k\in K} \dista_X(y, k),
\end{equation}
where $\dista_X$ is the metric of $X$.
Moreover, a set $\calA$ is the global (or ``universal'') 
attractor of the  semigroup $S(\cdot)$ iff $\calA$ is 
attracting and compact in $X$ and is fully invariant with
respect to $S(\cdot)$, i.e., $S(t)\calA = \calA$ for all $t\ge 0$.

For our system, the existence of the global attractor will be 
obtained by showing asymptotic compactness of the semigroup $S(\cdot)$, 
meaning, roughly speaking, that, for $t$ sufficiently large, all
solutions emanating from initial data satisfying \eqref{hp:init}
take values in a set that is bounded in some ``better space''.
To make precise this concept, also in view of the fact that 
the phase space $\Phi_m\times \Sigma$ and the distance 
\eqref{dist:Fhim} are somehow nonstandard, we introduce some more
tools, again following the lines of the approach devised in 
\cite{RS}. As before, we shall consider the behavior of $\fhi$ 
and $\sigma$ ``separately'', and, complementing respectively 
\eqref{defi:Fhim} and \eqref{defi:Sigma}, we define 
\begin{align}\label{defi:Psim}
  \Psi_m:=\big\{ \fhi \in W:~\beta^0(\fhi) \in L^2(\Omega),~\fhi\OO=m\big\},\\
 \label{defi:V+}
  V_+:=\big\{ \sigma\in V= D(A^{1/2}):~
   \sigma\ge 0~\text{a.e.~in}~\Omega\big\}.
\end{align}
We recall that, for $r\in \dom(\beta)$, $\beta^0(r)$ denotes the element 
of minimum modulus of the convex set $\beta(r)$ (recall that, according to
the general theory detailed in \cite{Ba,Br}, we may consider multivalued maximal 
monotone operators $\beta$).
\noindent%
Mimicking \eqref{dist:Fhim}, the set $\Psi_m$ can be naturally endowed 
with the distance 
\begin{equation}\label{dist:Psim}
  \dist_1(\fhi_1,\fhi_2) 
   := \| \fhi_1 - \fhi_2 \|_{H^2(\Omega)}
    + \| \beta^0(\fhi_1) - \beta^0(\fhi_2) \|.
\end{equation}
\beos\label{oss:RS}
As observed in \cite[Prop.~3.15]{RSS} (and it is also easy to check
directly), in the case when $\dom\betaciapo=[-1,1]$ 
(as happens for \eqref{Flog}), the second summand in \eqref{dist:Psim} 
can be avoided, and one can simply use the $H^2$-norm. 
Recall that, in that case, the structure of the 
phase space is simpler as well (cf.~Remark~\ref{rem:casolog}).
\eddos
\noindent%
Then, an easy adaptation of the argument in \cite[Prop.~3.15]{RS}
gives
\bepr\label{pro:RS}
 It holds that $\Psi_m \subset \Phi_m$ with compact immersion;
 namely, if $\{\fhi_n\}$ is a sequence in $\Psi_m$ bounded with respect
 to the distance \eqref{dist:Psim}, then there exist $\fhi\in\Phi_m$
 and a nonrelabelled subsequence of $\{\fhi_n\}$ such that 
 $\dist(\fhi_n,\fhi)\to 0$.
\empr
\noindent%
We can then prove the following result providing existence of 
the global attractor in the regularity setting of 
Theorem~\ref{teo:wellpo}:
\bete\label{teo:attra}
 Let\/ assumptions {\rm (A1)-(A2)} hold together with the compatibility
 condition\/ \eqref{keycoeff}. Then, the dynamical process $S(\cdot)$
 associated to system~\eqref{CH1}-\eqref{nutr} admits the global attractor
 $\calA\subset \Psi_m\times V_+$. More precisely,  $\calA$ is bounded 
 in $\Psi_m\times V_+$, in the sense that there exists $M>0$ such that 
 \begin{equation}\label{bound:A}
   \dist_1(\phi,0) + \| s \|_V \le M,
    \quext{for every }\,(\phi,s)\in \calA.
 \end{equation}
 Moreover, for what concerns the phase variable, one has 
 more precisely
 \begin{equation}\label{bound:A1}
   \| \phi \|_{W^{2,6}(\Omega)} 
    + \| \beta^0(\phi) \|_{L^{6}(\Omega)} \le M,
    \quext{for every }\,(\phi,s)\in \calA.
 \end{equation}
\ente
\noindent%
Exhibiting stronger a priori estimates, we can also prove some additional 
regularity of the nutrient states in the attractor. Namely, we have
\bete\label{teo:regattra}
 Let\/ the assumptions of Theorem~\ref{teo:attra} hold.
 Then, the elements of the attractor $\calA$ satisfy the following 
 additional regularity estimate: there exists $M_1>0$ such that 
 \begin{equation}\label{bound:A2}
    \| s \|_{H^2(\Omega)} \le M_1,
        \quext{for every }\,(\phi,s)\in \calA.
 \end{equation}
\ente
\noindent%
\beos\label{oss:regoattr}
It is worth observing that the regularity properties of $\fhi$ 
stated in \eqref{bound:A1} correspond to what is currently 
considered the ``state of the art'' for the Cahn-Hilliard equation 
with singular potential, at least in three dimensions of space
(cf., e.g., \cite{MiZe}). On the other hand,
regarding $\sigma$, in the ``high-regularity regime'' \eqref{nutr}
behaves more or less like a semilinear parabolic equation;
hence the only real obstacle that prevents us from further improving 
the asymptotic regularity of $\sigma$ comes from the forcing 
terms depending on $\fhi$, whose smoothness, as said,
is limited by \eqref{bound:A1}.
\eddos
\noindent%
The last question we would like to discuss is related to the 
nonnegativity of $\sigma$, which is incorporated in the definition
\eqref{defi:Sigma} of the phase space. Actually, as noted in 
Remark~\ref{rem:posi}, using the Stampacchia truncation method
one can easily prove that such a property is maintained in time.
On the other hand, it may be worth investigating whether this 
``minimum principle'' could be improved, possibly in an uniform way 
with respect to time, by means of parabolic smoothing arguments.
Actually, as observed in some previous papers
dealing with related models (see, e.g., \cite{LRS,Ssens}),
it is very easy to prove that, if one additionally assumes
\begin{equation}\label{strongpos}
  \sigma_0>0\quext{a.e.~in~}\, \Omega,
   \qquad \ln \sigma_0\in L^1(\Omega),
\end{equation}
then also this property is conserved in time. We also remark that
this is far from being just an ``academic'' question, as the regularity 
\eqref{strongpos} may serve as a source of additional a-priori estimates,
on which is based, for instance, the concept of ``entropic'' solution 
developed in~\cite{LRS}.

Here, assuming \eqref{strongpos}, we can show a time 
regularization property for $\ln \sigma$ implying that
for strictly positive times $\sigma(t)$ is separated from $0$
in the uniform norm. More precisely, we have
\bete\label{teo:ln}
 Let the assumptions of\/ {\rm Theorem~\ref{teo:wellpo}} hold and let,
 in addition,\/ \eqref{strongpos} be satisfied.
 %
 %
 Then, for any $0<\tau<T<\infty$ there exists $\delta=\delta(\tau,T)>0$
 such that 
 \begin{equation}\label{strongpos3} 
   \sigma(x,t) \ge \delta 
    \quext{for every }\,(x,t)\in \Omega\times[\tau,T].
 \end{equation}
\ente
\noindent%
\beos\label{oss:nodiss}
 As noticed in the introduction, it is not obvious whether an estimate 
 of the form \eqref{strongpos} could be proved uniformly in time. It is actually
 clear that, if the initial datum satisfies $\sigma_0(x) \ge \delta_0$ 
 a.e.~in~$\Omega$ for some constant $\delta_0>0$, then \eqref{strongpos}
 holds with $\tau=0$. On the other hand, the a-priori 
 estimates leading to that relation apparently lack a dissipative character,
 whence one may expect that $\delta$ in \eqref{strongpos} could degenerate
 to $0$ as one lets $T\nearrow \infty$. As a consequence it is also not
 clear whether the states $\sigma$ in the attractor are separated
 from $0$ in the uniform norm.
\eddos


Finally, for later convenience, it is useful to recall the so-called {\sl uniform Gr\"onwall
lemma}\/ (cf., e.g., \cite[Lemma~III.1.1]{Te}):
\bele\label{lem:uni}
 Let $y,a,b\in L^1_{\loc}(0, \infty)$ be three non-negative functions 
 such that 
 \begin{equation}\label{uni:gro}
   y'(t) \le a(t) y(t) + b(t) \quext{for a.e.\ }\, t > 0, 
 \end{equation}  
 and let $a_1$, $a_2$, $a_3$ three non-negative constants such that
 \begin{equation}\label{uni:gro2}
   \sup_{t\ge0} \int_t^{t+1} a(s) \le a_1, \qquad 
    \sup_{t\ge0} \int_t^{t+1} b(s) \le a_2, \qquad 
    \sup_{t\ge0} \int_t^{t+1} y(s) \le a_3.
 \end{equation}  
 Then, it follows that $y(t + 1) \le (a_2 + a_3)e^{a_1}$
 for every $t > 0$.
\enle


\section{A priori estimates, existence and dissipativity}
\label{sec:apriori}

In this part we detail the basic a priori estimates leading to existence 
of weak solutions in the regularity setting of Theorem~\ref{teo:wellpo}.
Our argument partially follows the lines of the proofs given in 
previous contributions on related models, see, e.g., 
\cite{GHW,RSS,Ssens}. For this reason we will limit 
ourselves to detail the basic a priori estimates needed for proving 
existence, by formulating them in a way that could also yield
the dissipativity of the associated dynamical process, which is 
a fundamental step of our long-time analysis. Moreover, we will work
directly on the ``original'' system \eqref{CH1}-\eqref{nutr} 
without referring to any regularization scheme. Actually, the
problem of the approximation 
(or discretization) of the system has been thoroughly discussed in 
the quoted references. However, we would like at least to recall a basic
point: for this class of systems, the coercivity of the energy
functional $\calE$ (cf.~\eqref{energy} below) is linked to the 
choice of a ``singular'' potential (as in Assumption~(A2)),
which guarantees a-priori uniform boundedness of $\fhi$ and, consequently, 
a control of the coupling term $-\chi\sigma\fhi$. Any approximation
of the system must take this fact into account: for instance, if
$\beta$ is replaced by a ``smooth'' regularization, as is customary
for Cahn-Hilliard models, then the coercivity
of $\calE$ is lost unless some further regularizing term is added 
in order to compensate the loss of boundedness of $\fhi$.
Actually, in our view, this seems to be, the main difficulty 
arising when one approximates models in this class.


\subsection{Dissipative energy estimate}
\label{subsec:en}

We derive here a variant of the energy estimate (in itself, this can 
be seen as a direct consequence of the variational structure of the model),
designed so to provide as a byproduct the existence of a bounded 
absorbing set with respect to the product metric in 
$\Phi_m\times \Sigma$. To obtain it, we proceed through several
steps. First of all, we observe the mass conservation property
\begin{equation}\label{st:mass}
  \ddt \io \fhi = 0,
\end{equation}
obtained integrating \eqref{CH1} over $\Omega$ and exploiting the no-flux
conditions. Next, we test \eqref{CH1} by $\mu$, \eqref{CH2} by $\fhi_t$, 
and \eqref{nutr} by $\ln \sigma - \chi \fhi$. Combining the resulting relations, 
it is not difficult to arrive at
\begin{equation}\label{st:11}
  \ddt \calE + \| \nabla \mu \|^2
   + \io \sigma | \nabla (\ln \sigma - \chi \fhi) |^2
   + \io \big( h(\sigma,\fhi) \sigma^2 
         - k(\sigma,\fhi) \sigma \big) (\ln \sigma - \chi\fhi) = 0,
\end{equation}
with the physical energy $\calE$ being given by
\begin{equation}\label{energy}
  \calE(\fhi,\sigma) = \frac12 \| \nabla\fhi \|^2 
   + \io \Big( F(\fhi) + \sigma (\ln \sigma - 1) - \chi \sigma\fhi \Big).
\end{equation}
Now, owing to condition \eqref{hk:sign} and to the fact 
\begin{equation}\label{fhibdd}
  \| \fhi \|_{L^\infty(0,T;L^\infty(\Omega))} \le 1,
\end{equation}
which, as noted before, is a direct consequence of Assumption~(A2)
at least as far as we work on the ``original'' system and not on a
regularization of it, it is not difficult to deduce 
\begin{equation}\label{st:12}
  \big( h(\sigma,\fhi) \sigma^2 - k(\sigma,\fhi) \sigma \big) (\ln \sigma - \chi\fhi) 
  \ge \frac{\hgiu}2 \sigma^2\ln (1 + \sigma) + \kappa \sigma \ln( 1 + \sigma) - c.
\end{equation}
Here and below, $c\ge 0$ and $\kappa>0$ are computable constants depending on the 
constants $\ksu,\kgiu,\hgiu$ and on the other assigned parameters of the system,
with $\kappa$ being used in estimates from below. The value of $c,\kappa$ may 
vary on occurrence. Of course, it is clear that, if the a priori estimates
were adapted to a (hypothetical) approximation, $c,\kappa$ would not be allowed
to depend on the approximation parameter(s).

That said, it is not difficult to check that the energy satisfies the coercivity 
property
\begin{equation}\label{coerc:en}
  \mathcal{E}( \varphi, \sigma) 
    \ge \kappa \big( \| \sigma \ln (1 + \sigma) \|_{L^1(\Omega)}
       + \| \fhi \|_{V}^2
       + \| \betaciapo(\fhi) \|_{L^1(\Omega)} \big)
       - c,
\end{equation}
as well as the boundedness condition
\begin{equation}\label{bound:en}
  \calE( \varphi, \sigma) 
    \le c \big( \| \sigma \ln (1 + \sigma) \|_{L^1(\Omega)}
       +  \| \fhi \|_{V}^2
       + \| \betaciapo(\fhi) \|_{L^1(\Omega)} + 1 \big).
\end{equation}
In order to obtain a dissipative version of the energy estimate,
we multiply \eqref{CH2} by $\fhi - \fhi\OO$, where we recall that 
$\fhi\OO$ stands for the (conserved, by \eqref{st:mass}) spatial 
average of $\fhi$. 
Then, 
%
%
let us recall that (cf.\ \cite{MiZe}) 
\begin{equation}\label{co:MZ}
  \io f(\fhi) (\fhi - \fhi\OO) 
   \ge \kappa \io \betaciapo(\fhi) + \kappa \io |\beta(\fhi)| - c,
\end{equation}
where $\kappa>0$ and $c\ge 0$ may now depend also on the ``assigned'' 
quantity $m=(\fhi_0)\OO$ and on the ``nonconvexity'' parameter $\lambda$.
Hence, we deduce 
\begin{equation}\label{co:21}
  \kappa \io \betaciapo(\fhi) + \kappa \io |\beta(\fhi)| 
   + \| \nabla \fhi \|^2 
  \le c 
   + \io \mu (\fhi - \fhi\OO)
   + \chi \io \sigma (\fhi - \fhi\OO).
\end{equation}
Now, using the Poincar\'e-Wirtinger inequality there follows
\begin{equation}\label{co:22}
  \io \mu (\fhi - \fhi\OO)
   = \io ( \mu - \mu\OO) (\fhi - \fhi\OO)
   \le \frac12 \| \nabla \fhi \|^2 
   + c_\Omega \| \nabla \mu \|^2,
\end{equation}
where $c\OO$ is an embedding constant. Then, replacing \eqref{co:22} into
\eqref{co:21}, multiplying the result by $1/(2c\OO)$ and adding it to 
\eqref{st:11}, using also \eqref{st:12} we obtain
\begin{align}\no
  & \ddt \calE + \frac12 \| \nabla \mu \|^2
   + \io \sigma | \nabla (\ln \sigma - \chi \fhi) |^2
   + \kappa \bigg( \io \betaciapo(\fhi) + \io |\beta(\fhi)| 
   + \| \nabla \fhi \|^2 \bigg)\\
 \label{co:24}
  & \mbox{}~~~~~
    + \frac{\hgiu}2 \io \sigma^2\ln (1 + \sigma) + \kappa \io \sigma \ln( 1 + \sigma) 
      \le c + \frac{\chi}{2c\OO} \io \sigma(\fhi - \fhi\OO),
\end{align}
for a possibly new value of the constant $\kappa>0$. 
To control the last term in \eqref{co:24}, we use \eqref{fhibdd}
so to obtain
\begin{equation}\label{co:23}
  \frac{\chi}{2c\OO} \io \sigma(\fhi - \fhi\OO)
   \le c \io \sigma 
   \le \frac{\hgiu}4 \io \sigma^2\ln (1 + \sigma) - c.
\end{equation}
We also observe that, as $|\fhi| \le 1$ almost everywhere,
there follows
\begin{equation}\label{fbeta}
  \io \betaciapo(\fhi) - \frac\lambda2|\Omega|
   \le \io F(\fhi) \le \io \betaciapo(\fhi).
\end{equation}
Hence, using $F$ or $\betaciapo$ in the estimates is in fact equivalent
(a similar equivalence holds for $f$ and $\beta$).
That said, neglecting the (positive) contribution of the 
cross-diffusion term and rearranging, \eqref{co:24} implies in particular
\begin{align}\no
  & \ddt \calE 
   + \kappa \bigg( \io \betaciapo(\fhi)  + \| \fhi \|_V^2 
   + \|  \sigma \ln( 1 + \sigma) \|_{L^1(\Omega)} + 1 \bigg)\\
 \label{co:24b}
  & \mbox{}~~~~~
   + \kappa \io |\beta(\fhi)|
   + \frac12 \| \nabla \mu \|^2
   + \frac{\hgiu}4 \| \sigma^2\ln (1 + \sigma) \|_{L^1(\Omega)}
      \le c,
\end{align}
or, in other words, owing to \eqref{bound:en},
\begin{equation}\label{co:25}
  \ddt \calE 
   + \kappa \calE 
   + \kappa_0 \| \nabla \fhi \|^2
   + \kappa \io |\beta(\fhi)|
   + \frac12 \| \nabla \mu \|^2
   + \frac{\hgiu}4 \| \sigma^2\ln (1 + \sigma) \|_{L^1(\Omega)}
      \le c,
\end{equation}
for possibly new values of the ``small'' positive constants $\kappa$ and 
$\kappa_0$ (we preferred to keep the latter as a separate value
as it will be managed in a while), which is a dissipative estimate 
for the physical energy $\calE$ and may be used for the construction
of the absorbing set. For this purpose, however, it will be convenient
to complement the above inequality with a further contribution.
Namely, we test \eqref{CH2} by $-\Delta\fhi$ so to obtain
\begin{equation}\label{co:33}
  \| \Delta \fhi \|^2 
   \le \lambda \| \nabla \fhi \|^2 
   + \| \nabla \mu \| \| \nabla \fhi \|
   + \chi \| \sigma \| \| \Delta \fhi \|,
\end{equation}
whence, by simple manipulations, 
\begin{equation}\label{co:33b}
  \| \Delta \fhi \|^2 
   \le (2 \lambda + 1) \| \nabla \fhi \|^2 
   + \| \nabla \mu \|^2
   + \chi^2 \| \sigma \|^2.
\end{equation}
Then, multiplying \eqref{co:33b} by $\delta>0$
taken sufficiently small and summing the result to \eqref{co:25},
it is not difficult to arrive at
\begin{equation}\label{co:25b}
  \ddt \calE 
   + \kappa \calE 
   + \frac{\kappa_0}2 \| \nabla \fhi \|^2
   + \kappa \io |\beta(\fhi)|
   + \kappa \| \Delta\fhi \|^2
   + \frac14 \| \nabla \mu \|^2 
   + \frac{\hgiu}8 \| \sigma^2\ln (1 + \sigma) \|_{L^1(\Omega)}
      \le c_0,
\end{equation}
where $c_0>0$ is a computable constant depending only on the
assigned parameters of the system and independent of the initial data.

Let us now assume that the initial data $(\fhi_0,\sigma_0)$ are chosen,
as in the statement of Theorem~\ref{teo:wellpo}, in a bounded
``ball'' $B={\ov{B}}(0,R)$ of the phase space of arbitrarily large, but
otherwise assigned, radius $R>0$, namely we have
\begin{equation}\label{initR}
  \dist(\fhi_0,0) + \| \sigma_0 \|_{D(A^{1/4})} \le R.
\end{equation}
By \eqref{bound:en}, the above implies in particular that there exists
a computable monotone function $Q:\RR^+ \to \RR^+$ such that 
$\calE(0)\le Q(R)$. Generally speaking, here and below
$Q$ will always denote a positive-valued function assumed to be
monotone with respect of each of its arguments. The expression 
of $Q$, may only depend on the fixed parameters of the system
and may vary on occurrence.

Applying the Gr\"onwall lemma to \eqref{co:25b}, we then deduce 
\begin{equation}\label{di:11}
  \calE(t) \le \calE(0) e^{-\kappa t} + \frac{c_0}{\kappa}
   \le Q(R) e^{-\kappa t} + \frac{c_0}{\kappa},
\end{equation}
for every $t\ge 0$. In addition to that, integrating 
\eqref{co:25b} over $[t,t+1]$ for a generic $t\ge 0$, 
we deduce 
\begin{align}\no
  & \int_t^{t+1} \bigg( 
   \frac{\kappa_0}2 \| \nabla \fhi \|^2
   + \kappa \io |\beta(\fhi)|
   + \kappa \| \Delta\fhi \|^2
   + \frac14 \| \nabla \mu \|^2
   + \frac{\hgiu}8 \| \sigma^2\ln (1 + \sigma) \|_{L^1(\Omega)} \bigg)\\
 \label{di:12}
   & \mbox{}~~~~~~~~~~
   \le Q(R) e^{-\kappa t} + C_0, 
\end{align}
still for every $t\ge 0$. Here and below, $C_0>0$ is a computable 
constant, varying on occurrence, whose value is independent of the initial 
data. For instance, in the above formula one can take 
$C_0 = c_0 ( 1 + \kappa^{-1} )$. With this notation at hand,
estimate \eqref{di:11}, thanks also to \eqref{coerc:en}, gives
\begin{equation}\label{di:11b}
  \| \fhi(t) \|_V^2 
   + \| \betaciapo(\fhi) \|_{L^1(\Omega)}
   + \| \sigma \ln (1 + \sigma) \|_{L^1(\Omega)}
   \le Q(R) e^{-\kappa t} + C_0.
\end{equation}
Note that the above is sufficient in order to control
$\fhi$ with respect to the ``energy distance''. 
On the other hand, regarding $\sigma$, the 
information provided by \eqref{di:11b} is very weak. 
Actually, as \eqref{defi:Sigma} suggests, the conditions on 
$\sigma$ required for having well-posedness
are much more restrictive compared to the information 
contained in the energy functional. Hence, in the sequel, 
we will need to derive some further {\sl dissipative}\/ estimates 
of better norms of $\sigma$. This procedure will be partially based on 
tools and techniques proposed in former papers dealing with 
related models.

%
%
%
%


\subsection{Estimate of the singular potential}
\label{subsec:log}

We go back to \eqref{co:21}, where we can now neglect some nonnegative terms 
on the \lhs, and modify a bit the control of the integral terms on the \rhs.
Using repeatedly \eqref{di:11b}, we first observe that
\begin{equation}\label{co:23c}
  \chi \io \sigma(\fhi - \fhi\OO)
   \le c \io \sigma 
   \le Q(R) e^{-\kappa t} + C_0.
\end{equation}
Similarly,
\begin{equation}\label{co:22b}
  \io \mu (\fhi - \fhi\OO)
   \le c \| \nabla \fhi \| \| \nabla \mu \|
   \le \big( Q(R) e^{-\kappa t} + C_0 \big)
        \| \nabla \mu \|.
\end{equation}
Moreover, using \eqref{co:23c}-\eqref{co:22b} and integrating \eqref{CH2}
over $\Omega$, it is not difficult to deduce from \eqref{co:21} that
\begin{equation}\label{mz:11}
  \| \beta(\fhi) \|_{L^1(\Omega)} 
  + | \mu\OO |
   \le c \big( 1 + \| \nabla\mu \| \big)
   \big( Q(R) e^{-\kappa t} + C_0 \big).
\end{equation}
%
%
%
Squaring, integrating over a generic interval $(t,t+1)$ and using 
the control of $\|\nabla\mu\|^2$ resulting from~\eqref{di:12}, 
we also obtain
\begin{equation}\label{mz:13}
  \int_t^{t+1} | \mu\OO |^2
   \le Q(R) e^{-\kappa t} + C_0,
\end{equation}
whence, using \eqref{di:12} once more, we deduce
\begin{equation}\label{di:14}
  \int_t^{t+1} \| \mu \|_V^2
   \le Q(R) e^{-\kappa t} + C_0.
\end{equation}
From the above, using \eqref{CH1} and standard arguments
based on \eqref{nor:*}-\eqref{nor:*2}, we also infer
\begin{equation}\label{di:15}
  \int_t^{t+1} \| \fhi_t \|_*^2
   \le Q(R) e^{-\kappa t} + C_0.
\end{equation}
%


\subsection{``Unilateral'' estimate of the singular term}
\label{subsec:unila}

We repeat here the argument devised in \cite[Sec.~4]{RSS}. Namely,
observing that, by \eqref{dombeta}, $\beta(r)$ has the same sign of $r$,
we test \eqref{CH2} by $-(\beta(\fhi)_-)^5$ (i.e., minus the fifth power
of the negative part of $\beta(\fhi)$). As this is a monotone function
of $\fhi$, the Laplacian, once integrated by parts, provides a nonnegative
contribution. Hence we end up with
\begin{equation}\label{mo:11}
  \| \beta(\fhi)_- \|_{L^6(\Omega)}^6
    = - \io \mu( \beta(\fhi)_-)^5
     - \io \lambda \fhi ( \beta(\fhi)_-)^5
    - \chi \io \sigma( \beta(\fhi)_-)^5.
\end{equation}
Now, the last term on the \rhs\ is obviously nonpositive, while the other two
integrals are controlled by observing that
\begin{equation}\label{mo:12}
   - \io \mu( \beta(\fhi)_-)^5 - \io \lambda \fhi ( \beta(\fhi)_-)^5
    \le c \big( \| \mu \|_{L^6(\Omega)} + \| \fhi \|_{L^6(\Omega)} \big)
    \| \beta(\fhi)_- \|_{L^6(\Omega)}^5.
\end{equation}
Hence, by the uniform boundedness of $\fhi$ and Sobolev's embeddings, 
we end up with 
\begin{equation}\label{mo:13}
  \| \beta(\fhi)_- \|_{L^6(\Omega)} 
    \le c \big( \| \mu \|_{V} + 1 \big),
\end{equation}
whence, squaring and integrating over $(t,t+1)$ for a generic $t\ge 0$, 
using \eqref{di:14}, we obtain
\begin{equation}\label{mo:14}
  \int_t^{t+1} \| \beta(\fhi)_- \|_{L^6(\Omega)}^2
   \le Q(R) e^{-\kappa t} + C_0.
\end{equation}
%


\subsection{Improved summability of $\sigma$}
\label{subsec:small}

In this part we adapt the argument devised in \cite{RSS} to 
obtain a ``decoupled'' estimate for the cross-diffusion term in 
\eqref{nutr}. This uses in an essential way the compatibility condition
\eqref{keycoeff}, which arises as the stabilizing contribution 
of the logistic source and the chemotactic effect provide terms
that have the same growth rate, so that the estimate
can be closed only by comparing the magnitude of the corresponding 
coefficients.

That said, we first observe that estimates \eqref{di:11}-\eqref{di:12} 
and \eqref{di:14} imply that, as the initial data satisfy 
condition \eqref{initR}, there exist a number $R_0$ independent of $R$ and a time 
$T_0$ depending on $R$, such that, for every $t\ge T_0$, there holds 
\begin{equation}\label{di:1112}
  \calE(t) + \int_t^{t+1} \bigg( 
   \| \fhi \|^2_V
   + \| \beta(\fhi) \|_{L^1(\Omega)}
   + \| \Delta\fhi \|^2
   + \| \mu \|_V^2
   + \| \sigma^2\ln (1 + \sigma) \|_{L^1(\Omega)} \bigg)
  \le R_0.
\end{equation}
%
%
%
Let us now test \eqref{nutr} by $\sigma^p$ for ``small'' $p>0$ to be chosen 
later on. Then, simple calculations give
\begin{equation}\label{di:22}
  \frac1{p+1} \ddt \| \sigma \|_{L^{p+1}(\Omega)}^{p+1}
   + \hgiu \| \sigma \|_{L^{p+2}(\Omega)}^{p+2}
   + \frac{4p}{(p+1)^2} \io | \nabla \sigma^{\frac{p+1}2} |^2
  \le \ksu \| \sigma \|_{L^{p+1}(\Omega)}^{p+1}
   + p \chi \io \sigma^p \nabla \sigma \cdot \nabla \fhi.
\end{equation}
Now, using \eqref{CH2}, the last term is managed as follows:
\begin{align}\no
  p \chi \io \sigma^p \nabla \sigma \cdot \nabla \fhi
   & = - \chi \frac{p}{p+1} \io \sigma^{p+1} \Delta \fhi\\
 \label{di:23} 
   & = \chi \frac{p}{p+1} \io \mu \sigma^{p+1}
   - \chi \frac{p}{p+1} \io f(\fhi) \sigma^{p+1}
   + \chi^2 \frac{p}{p+1} \io \sigma^{p+2}.
\end{align}
Following \cite{RSS}, the first two terms on the \rhs\ 
can be controlled by observing, respectively, that 
\begin{equation}\label{di:25}
  \chi \frac{p}{p+1} \io \mu \sigma^{p+1}
   \le c \| \mu \|_{L^4(\Omega)} \| \sigma^{\frac{p+1}2} \| 
            \| \sigma^{\frac{p+1}2} \|_{L^4(\Omega)}
   \le \frac{p}{(p+1)^2} \| \sigma^{\frac{p+1}2} \|_V^2
    + c_p \| \mu \|_{V}^2 \| \sigma^{\frac{p+1}2} \|^2
\end{equation}
and that (here we also use that $\sigma$ is nonnegative)
\begin{align}\no
  - \chi \frac{p}{p+1} \io f(\fhi) \sigma^{p+1}
   & = - \chi \frac{p}{p+1} \io \beta(\fhi) \sigma^{p+1}
    + \lambda \chi \frac{p}{p+1} \io \fhi \sigma^{p+1}\\
 \no    
   & \le  \| \beta(\fhi)_- \|_{L^6(\Omega)} \| \sigma^{\frac{p+1}2} \|
       \| \sigma^{\frac{p+1}2} \|_{L^3(\Omega)}
     + c \| \sigma \|_{L^{p+1}(\Omega)}^{p+1}\\
 \label{di:26}
   & \le \frac{p}{(p+1)^2} \| \sigma^{\frac{p+1}2} \|_V^2
    + c_p \| \beta(\fhi)_- \|_{L^6(\Omega)}^2 \| \sigma^{\frac{p+1}2} \|^2
    + c \| \sigma \|_{L^{p+1}(\Omega)}^{p+1}.
\end{align}
Collecting \eqref{di:23}-\eqref{di:26}, \eqref{di:22} gives
\begin{align}\no
  & \frac1{p+1} \ddt \| \sigma \|_{L^{p+1}(\Omega)}^{p+1}
   + \Big[ \hgiu - \chi^2 \frac{p}{p+1} \Big] \| \sigma \|_{L^{p+2}(\Omega)}^{p+2}
   + \frac{2p}{(p+1)^2} \io | \nabla \sigma^{\frac{p+1}2} |^2\\
 \label{di:22b}  
  & \mbox{}~~~~~ 
  \le c_p \Big( 1 + \| \mu \|_{V}^2 + \| \beta(\fhi)_- \|_{L^6(\Omega)}^2 \Big)
   \| \sigma \|_{L^{p+1}(\Omega)}^{p+1}.
\end{align}
The above procedure works for every $p\in(0,1]$, but
provides an additional estimate only if 
\begin{equation}\label{hchi}
  \hgiu - \chi^2 \frac{p}{p+1} \ge 0.
\end{equation}
Actually, in the paper \cite{RSS} it was directly assumed that 
the above condition was verified for $p=1$, which corresponds to the 
fact that the chemotactic response coefficient $\chi$ is ``small''
compared to the logistic source.
Here, we may observe that, as $\hgiu$ is assigned, it is clear that 
\eqref{hchi} holds at least for $p>0$ sufficiently small. On the other 
hand, \eqref{di:22b} seems to provide some regularity gain only for
$p\ge 1/2$, which corresponds exactly to condition \eqref{keycoeff}.

In the sequel, we will directly assume the worst scenario
we are able to deal with, i.e.\ the case when \eqref{hchi}
holds exactly for $p=1/2$ (or, in other words, equality holds 
in \eqref{keycoeff}), noting that the argument can be simplified
if one directly assumes that \eqref{hchi} holds for $p=1$ as
done in \cite{RSS}. For $p=1/2$, 
using \eqref{di:22b} and the uniform 
Gr\"onwall lemma (Lemma~\ref{lem:uni}), thanks also to \eqref{di:14}
and \eqref{mo:14}, for every $t\ge T_0 + 1$ (where $T_0$
was defined at the beginning of this subsection) one deduces that 
\begin{equation}\label{di:41}
  \| \sigma(t) \|_{L^{3/2}(\Omega)}^{3/2}
   + \int_t^{t+1} |\nabla \sigma^{3/4} |^2
   \le Q(R_0),
\end{equation}
which, by Sobolev's embeddings and elementary interpolation, 
implies 
\begin{equation}\label{di:42}
  \| \sigma \|_{L^2(t,t+1;L^3(\Omega))}
   \le Q(R_0).
\end{equation}


\subsection{``Bilateral'' estimate of the singular term}
\label{subsec:bila}

We proceed similarly with Subsec.~\ref{subsec:unila}, namely we test 
\eqref{CH2} by $|\beta(\fhi)|\beta(\fhi)$, now with no need for
taking the negative part. Actually, 
the choice of a quadratic function is motivated 
by the $L^3$- regularity in \eqref{di:42}. Using the monotonicity
of the function $r\mapsto | \beta(r) | \beta(r)$, we have
\begin{equation}\label{mo:21}
  \| \beta(\fhi) \|_{L^3(\Omega)}^3
    \le \io \mu | \beta(\fhi) | \beta(\fhi)
    + \io \lambda \fhi | \beta(\fhi) | \beta(\fhi)
    + \chi \io \sigma | \beta(\fhi) | \beta(\fhi).
\end{equation}
Now, proceeding similarly with Subsec.~\ref{subsec:unila}, it is 
not difficult to deduce 
\begin{equation}\label{mo:22}
  \| \beta(\fhi) \|_{L^3(\Omega)} 
   \le c\big( 1 + \| \mu \|_{L^3(\Omega)} + \| \sigma \|_{L^3(\Omega)} \big).
\end{equation}
Then, squaring the above relation and performing a further comparison of 
terms in \eqref{CH2}, we get
\begin{equation}\label{mo:23}
  \| \beta(\fhi) \|_{L^3(\Omega)}^2
  + \| \Delta\fhi \|_{L^3(\Omega)}^2
   \le c\big( 1 + \| \mu \|^2_{L^3(\Omega)} + \| \sigma \|^2_{L^3(\Omega)} \big).
\end{equation}
Hence, integrating over $(t,t+1)$ and using \eqref{di:14}, \eqref{di:42},
and elliptic regularity results, we deduce 
\begin{equation}\label{mo:24}
  \| \beta(\fhi) \|_{L^2(t,t+1;L^3(\Omega))}
   + \| \fhi \|_{L^2(t,t+1;W^{2,3}(\Omega))}
   \le Q(R_0),  
\end{equation}
for every $t\ge T_0 + 1$.
\beos\label{rem:L6}
If \eqref{hchi} holds for $p=1$, it is easy to see that,
modifying the above argument, one can directly get $L^6$-
and $W^{2,6}$- space summability in \eqref{mo:23}. We will arrive 
at the same conclusion also for $p=1/2$, but in this case the proof is
a bit more involved.
\eddos


\subsection{``Parabolic'' regularity of $\sigma$}
\label{subsec:para}

We now repeat \eqref{di:22} with $p=1$ and integrate by parts the last term
as in \eqref{di:23}. With \eqref{mo:23} at hand, then one 
simply observes that
\begin{equation}\label{di:23b}
  - \frac{\chi}{2} \io \sigma^2 \Delta \fhi
   \le \| \sigma \| \| \sigma \|_{L^6(\Omega)} \| \Delta \fhi \|_{L^3(\Omega)} 
   \le \frac14 \| \sigma \|_{V}^2 
    + c \| \Delta \fhi \|_{L^3(\Omega)}^2 \| \sigma \|^2.
\end{equation}
Hence, \eqref{di:22}, for $p=1$, gives
\begin{equation}\label{di:22c}
  \frac1{2} \ddt \| \sigma \|^2
   + \hgiu \| \sigma \|_{L^{3}(\Omega)}^3
   + \frac{3}4 \io | \nabla \sigma |^2
  \le c \big( 1 + \| \Delta \fhi \|_{L^3(\Omega)}^2\big)
   \| \sigma \|^2.
\end{equation}
Then, 
applying once more the uniform 
Gr\"onwall lemma, using also \eqref{mo:23}, we deduce 
\begin{equation}\label{mo:25}
  \| \sigma(t) \|
   + \| \sigma \|_{L^3(t,t+1;L^3(\Omega))}
   + \| \sigma \|_{L^2(t,t+1;V)}
   \le Q(R_0),  
\end{equation}
say, for every $t\ge T_0+2$.

The above condition, in turn, permits us to further improve the 
regularity of $\fhi$ as anticipated in Remark~\ref{rem:L6}.
To this aim, we proceed similarly with 
Subsec.~\ref{subsec:bila}, now with an improved exponent.
Namely, we test \eqref{CH2} by $|\beta(\fhi)|^4\beta(\fhi)$
so to obtain the following analogue of \eqref{mo:23}:
\begin{equation}\label{mo:26}
  \| \beta(\fhi) \|_{L^6(\Omega)}^2
  + \| \Delta\fhi \|_{L^6(\Omega)}^2
   \le c\big( 1 + \| \mu \|^2_{L^6(\Omega)} + \| \sigma \|^2_{L^6(\Omega)} \big).
\end{equation}
Hence, integrating once more over $(t,t+1)$ and using now 
\eqref{di:14}, \eqref{mo:25} and elliptic regularity results, 
we arrive at 
\begin{equation}\label{mo:27}
  \| \beta(\fhi) \|_{L^2(t,t+1;L^6(\Omega))}
   + \| \fhi \|_{L^2(t,t+1;W^{2,6}(\Omega))}
   \le Q(R_0),  
\end{equation}
at least for $t\ge T_0+2$.


\subsection{Fractional parabolic regularity for $\sigma$}
\label{subsec:frac}

Using the operator $A$ introduced in \eqref{defi:A}, equation \eqref{nutr} can
be rewritten as 
\begin{equation}\label{nutr:A}
  \sigma_t + A \sigma + \chi \dive (\sigma\nabla\fhi) =
   - h(\sigma,\fhi) \sigma^2 + ( k(\sigma,\fhi) + 1) \sigma.
\end{equation}
Then, setting $H:= - h(\sigma,\fhi) \sigma^2 + ( k(\sigma,\fhi) + 1) \sigma$,
using Assumption~(A1), \eqref{mo:25} and interpolation, it is easy to check that
\begin{equation}\label{para:11}
   \| H \|_{L^2(t,t+1;D(A^{-1/4}))}
    \le c \| H \|_{L^2(t,t+1;L^{3/2}(\Omega))}
    \le Q(R_0),
    \quext{for a.e.\ }\,t\ge T_0+2,
\end{equation}
where we used the dual embedding of
\begin{equation}\label{para:12}
  D(A^{1/4})= H^{1/2}(\Omega) \subset L^3(\Omega).
\end{equation}
Condition \eqref{para:11} allows us to multiply \eqref{nutr:A} 
by $A^{1/2}\sigma$ so to deduce 
\begin{equation}\label{para:13}
  \frac12 \ddt \| A^{1/4} \sigma \|^2
   + \| A^{3/4} \sigma \|^2
  = ( H, A^{1/2} \sigma ) 
   - \chi \io \nabla \sigma \cdot \nabla\fhi A^{1/2}\sigma
   - \chi \io \sigma \Delta\fhi A^{1/2} \sigma,
\end{equation}
and we need to control the terms on the \rhs. We first notice that 
\begin{equation}\label{para:14}
  ( H, A^{1/2} \sigma ) 
   \le \frac14 \| A^{3/4} \sigma \|^2 
    + c \| A^{-1/4} H \|^2.
\end{equation}
Next,
\begin{align}\no
   - \chi \io \nabla \sigma \cdot \nabla\fhi A^{1/2}\sigma
   & \le c \| \nabla \sigma \| \| \nabla \fhi \|_{L^\infty(\Omega)} 
         \| A^{1/2} \sigma \|\\
 \no
   & \le c \| \nabla \fhi \|_{L^\infty(\Omega)} \| A^{1/2} \sigma \|^2\\
 \no
   & \le c \| \nabla \fhi \|_{L^\infty(\Omega)} \| A^{1/4} \sigma \| \| A^{3/4} \sigma \|\\
 \label{para:15}
   & \le c \| \fhi \|_{W^{2,6}(\Omega)}^2 \| A^{1/4} \sigma \|^2 
    + \frac14 \| A^{3/4} \sigma \|^2.
\end{align}
Finally,
\begin{align}\no
   - \chi \io \sigma \Delta\fhi A^{1/2} \sigma
   & \le c \| \sigma \|_{L^3(\Omega)}  \| \Delta \fhi \|_{L^3(\Omega)} 
         \| A^{1/2} \sigma \|_{L^3(\Omega)}\\
 \no
  & \le c \| A^{1/4} \sigma \| \| \Delta \fhi \|_{L^3(\Omega)}
        \| A^{3/4} \sigma \|\\
 \label{para:16}
  & \le c \| A^{1/4} \sigma \|^2 \| \Delta \fhi \|_{L^3(\Omega)}^2
        + \frac14 \| A^{3/4} \sigma \|^2.
\end{align}
Collecting \eqref{para:14}-\eqref{para:16}, \eqref{para:13} gives
\begin{equation}\label{para:13b}
  \frac12 \ddt \| A^{1/4} \sigma \|^2
   + \frac14 \| A^{3/4} \sigma \|^2
  \le c \| A^{-1/4} H \|^2
   +  c \| \fhi \|_{W^{2,6}(\Omega)}^2 \| A^{1/4} \sigma \|^2.
\end{equation}
Hence, recalling \eqref{mo:25}, \eqref{mo:27} and 
\eqref{para:11}, a further application of the uniform Gr\"onwall 
lemma yields
\begin{equation}\label{mo:31}
  \| \sigma(t) \|_{D(A^{1/4})}
   + \| \sigma \|_{L^2(t,t+1;D(A^{3/4}))}
   \le Q(R_0),  
\end{equation}
at least for every $t\ge T_0 + 3$.

\smallskip
\noindent%
Estimate \eqref{mo:31} is the last property we need in order to conclude 
the existence and dissipativity parts of the proof of Theorem~\ref{teo:wellpo}.
Indeed, from \eqref{di:11b} and \eqref{mo:31}
it is clear that property \eqref{asso:2}
holds, where $R_1$ (i.e., the radius of the absorbing set
w.r.t.\ the product distance of $\Fhi_m\times D(A^{1/4}))$
is a computable function of the quantities occurring
on the \rhs s of  \eqref{di:11b} and \eqref{mo:31};
moreover, one can take, for instance, $T_1 :=T_0 + 3$ in 
accordance with the validity range of \eqref{mo:31}.
\beos\label{smalltimes}
As our focus is mostly on the long-time behavior of 
solutions, we decided to present the a-priori bounds in a way 
that emphasizes their dissipative character. In particular,
we may check that, in fact, the values $R_1$
and $T_1$ depend only on the norm 
$\| \sigma_0 \ln(1 + \sigma_0)\|_{L^1(\Omega)}$ and 
not on the (stronger) $D(A^{1/4})$- norm of $\sigma_0$.
On the other hand, in order to apply the 
theory of infinite-dimensional dynamical systems,
we also need to prevent (possible) nonuniqueness phenomena
occurring ``earlier'' than the time $T_1$.
For this reason, we also need to achieve the regularity 
corresponding (for instance) to \eqref{mo:31}
(but the same applies to all the intermediate estimates 
that have been obtained through the uniform Gr\"onwall
lemma) over the time interval $(0,T_1)$. Of course this
can be easily done by adapting the previous arguments
and using the standard (rather than the uniform) Gr\"onwall
lemma for small values of the time variable. In this way,
one can complement (for instance) \eqref{mo:31}
by the corresponding relation 
\begin{equation}\label{mo:31small}
  \| \sigma \|_{L^\infty(0,T_1;D(A^{1/4}))}
   + \| \sigma \|_{L^2(0,T_1;D(A^{3/4}))}
   \le Q(R,T_1),  
\end{equation}
where, now, the computable function $Q$ depends also on the 
$D(A^{1/4})$- norm of $\sigma_0$, 
and not only on $\| \sigma_0 \ln(1 + \sigma_0)\|_{L^1(\Omega)}$.
\eddos


\subsection{End of proof of existence and dissipativity}
\label{subsec:exi}

In this part we detail how the estimates obtained so far imply 
existence of solutions and dissipativity of the associated semigroup
in the regularity setting of Theorem~\ref{teo:wellpo}. 

First of all, we check \eqref{rego:fhi}-\eqref{rego:sigma}, and, 
as we have decided to avoid detailing an approximation
scheme, we will just see how these regularity properties follow
from the a-priori estimates. In particular, as observed in 
Remark~\ref{smalltimes}, we will refer to the ``dissipative''
estimates detailed before, but of course, working on a generic 
interval $(0,T)$, in fact we also need their ``nondissipative'' counterparts
holding for small values of time.

That said, the first and the last regularity conditions in 
\eqref{rego:fhi} follow, respectively, from \eqref{di:15},
\eqref{mo:27}, and their non-dissipative counterparts.
To get the third condition, we test \eqref{CH2} by 
$-\Delta\fhi$ so to obtain 
\begin{equation}\label{co:d1}
  \| \Delta\fhi\|^2 
   \le \| \nabla\fhi \| \big( \| \nabla \mu \| 
    + \chi \| \nabla\sigma \| + \lambda \| \nabla \fhi \| \big),
\end{equation}
whence, observing that $\| \nabla \fhi \|$ is controlled uniformly in
time by \eqref{di:11}, the desired control follows by squaring \eqref{co:d1}
and integrating in time. Finally, we note that \eqref{di:11} implies 
at least that $\fhi\in L^\infty(0,T;V)$. In addition to that,
\eqref{mo:27} entails the second of \eqref{rego:ffhi}. 

The conditions proved so far are sufficient in order to apply the analogue
of \cite[Lemma~4.1]{RS} (which in turn is based on the abstract result
\cite[Lemma~3.9]{RS}), which implies the continuity properties appearing 
in the second \eqref{rego:fhi} and in the first \eqref{rego:ffhi}
(in fact, the lemma implies the absolute continuity of the function
$t\mapsto\calE(\fhi,\sigma)$).

Next, \eqref{rego:mu} follows directly
from \eqref{di:14}. Finally, concerning \eqref{rego:sigma}, one 
can first of all take advantage of \eqref{mo:31} (and of the corresponding 
nondissipative relation \eqref{mo:31small}). Then, to obtain
the regularity of the time derivative, this follows simply 
by testing \eqref{nutr:A} by $A^{-1/2}\sigma_t$ and proceeding 
similarly with Subsec.~\ref{subsec:frac} in order to control the coupling
and logistic terms (we leave the details to the reader). Note 
that, as a byproduct of this procedure, we deduce also the 
the continuity property in \eqref{rego:sigma} (which corresponds
to the continuity of the $\sigma$-component of the trajectory
with respect to the topology of $\Sigma$). Indeed, we have 
\begin{equation}\label{co:d2}
   \frac12 \ddt \| A^{1/4} \sigma \|^2
    = ( A^{1/4} \sigma_t, A^{1/4} \sigma )
    = ( A^{-1/4} \sigma_t, A^{3/4} \sigma ),
\end{equation}
and we may observe that the \rhs\ is summable as a consequence of 
the regularity properties proved so far.

As all the regularity properties in the statement have been achieved,
it is easy to check that, at least formally, 
these suffice in order to interpret 
the system in the form \eqref{CH2} plus \eqref{neum:fhi}-\eqref{init}
(in particular, the initial data are recovered as both
$\fhi$ and $\sigma$ are continuous in time with respect to
suitable topologies). Finally, concerning dissipativity,
for what concerns $\fhi$ this is a direct consequence of 
\eqref{di:1112}, whereas for $\sigma$ we can use \eqref{mo:31}.
In particular, one may take $T_1=T_0 + 3$ and 
choose $R_1=Q(R_0)$ for some explicitly
computable expression of $Q$ also depending on the 
coercivity constant $\kappa>0$ in \eqref{coerc:en} and 
on the expression of the function $Q$ appearing in \eqref{mo:31}.
This completes the proof of Theorem~\ref{teo:wellpo} for what 
concerns existence and dissipativity. Of course, as already noted,
this is a simplified proof in the sense that we worked directly
on the ``original'' system; if an approximation scheme is 
used, then the procedure should be adapted in order to show 
that the estimates are sufficient to pass to 
the limit in the approximation. This seems, however, a very
standard procedure, also in view of the fact that we are 
considering here a rather high regularity setting.
It is worth observing that recovering existence of solutions
by passing to the limit in an approximation might be a more delicate 
task in a lower regularity setting, see, e.g., \cite{LRS}.


\section{Uniqueness}
\label{sec:uni}

We prove here the uniqueness part of Theorem~\ref{teo:wellpo}. The argument 
partially follows the lines of \cite[Proof of Theorem~2.8]{RSS}; however, 
as the current setting is different in several aspects
(lack of a source term in \eqref{CH1},
different hypotheses on the logistic terms in \eqref{nutr}, fractional
regularity setting for $\sigma$), we prefer to report the proof in full detail.
That said, we consider a couple of solutions $(\fhi_1,\mu_1,\sigma_1)$ 
and $(\fhi_2,\mu_2,\sigma_2)$, both emanating from the same initial 
datum $(\fhi_0,\sigma_0)$ satisfying \eqref{hp:init}, and we 
look for a contractive estimate. To obtain it, we first set 
$(\fhi,\mu,\sigma):=(\fhi_1,\mu_1,\sigma_1)-(\fhi_2,\mu_2,\sigma_2)$,
and, taking the difference of system \eqref{CH1}-\eqref{nutr},
we infer
\begin{align}\label{CH1:d}
  & \fhi_t - \Delta \mu = 0,\\
 \label{CH2:d}
  & \mu = - \Delta \fhi + \beta(\fhi_1) - \beta (\fhi_2) - \lambda \fhi - \chi \sigma,\\
 \label{nutr:d}
  & \sigma_t - \Delta\sigma + \chi \dive( \sigma_1 \nabla \fhi ) 
   + \chi \dive (\sigma \nabla \fhi_2)
   = \big( - h_1 \sigma_1^2 + h_2 \sigma_2^2 
   + k_1 \sigma_1 - k_2 \sigma_2 \big),
\end{align} 
where we have set $h_1:=h(\sigma_1,\fhi_1)$, $h_2:=h(\sigma_2,\fhi_2)$, 
and so on. We will work on a generic, but otherwise fixed, time interval
$(0,T)$ and assume that both solutions turn out to satisfy the regularity 
conditions \eqref{rego:fhi}-\eqref{rego:sigma} over $(0,T)$.

Then, observing that $\fhi\OO=0$ and recalling that $\calN$ denotes
the inverse Neumann Laplacian over the functions with zero spatial mean, 
we may test \eqref{CH1:d} by $\calN \fhi$ and \eqref{CH2:d} 
by $\fhi$. Taking the difference of the resulting relations,
and using \eqref{nor:*} together with the monotonicity of $\beta$,
we deduce
\begin{equation}\label{co:51}
  \frac12 \ddt \| \fhi \|_*^2
   + \| \nabla \fhi \|^2
   \le \lambda \| \fhi \|^2 + \chi \io \sigma \fhi
   = \lambda \| \fhi \|^2 + \chi \io (\sigma - \sigma\OO) \fhi,
\end{equation}
the last equality holding as $\fhi$ has zero spatial average. Then,
adding $\| \fhi \|^2$ to both hand sides so to recover the full
$V$-norm of $\fhi$ on the \lhs, we may observe that the 
\rhs\ can then be managed this way:
\begin{equation}\label{co:51b}
  (\lambda + 1) \| \fhi \|^2 + \chi \io (\sigma - \sigma\OO) \fhi
   \le \frac 12 \| \fhi \|_V^2 + c \| \fhi \|^2_*
    + c \| \sigma - \sigma\OO \|_*^2,
\end{equation}
where we have also used Ehrling's inequality. Hence, 
\eqref{co:51} gives
\begin{equation}\label{co:52}
  \ddt \| \fhi \|_*^2
   + \| \fhi \|^2_V
   \le c \| \fhi \|_*^2 + c \| \sigma - \sigma\OO \|_*^2.
\end{equation}
Dealing with equation \eqref{nutr:d} is a bit trickier as the spatial
mean of $\sigma$ is not conserved. To this aim,
we first estimate the difference
of the mean values by integrating \eqref{nutr:d} over $\Omega$
so to obtain
\begin{equation}\label{co:53}
  \ddt \sigma\OO
  = \frac1{|\Omega|} \io \big( - h_1 \sigma_1^2 + h_2 \sigma_2^2 
   + k_1 \sigma_1 - k_2 \sigma_2 \big).
\end{equation}
Now, thanks to the global Lipschitz continuity of $h$ (cf.\ \eqref{hk:reg}),
it is not difficult to see that 
\begin{align}\no
  \bigg| \io \big( - h_1 \sigma_1^2 + h_2 \sigma_2^2 \big) \bigg| 
   & = \bigg| - \io (h_1 - h_2) \sigma_1^2 
   - \io h_2 (\sigma_1 + \sigma_2) \sigma \bigg| \\
 \no
   & \le c \| \sigma_1 \|_{L^3(\Omega)} \| \sigma_1 \|_{L^6(\Omega)} 
     \big( \| \sigma \| + \| \fhi \| \big)
   + c \big( \|\sigma_1 \| + \| \sigma_2 \| \big) \| \sigma \|,\\
 \label{co:54}
   & \le C \big( \| \sigma_1 \|_{V} + 1 \big)
     \big( \| \sigma \| + \| \fhi \| \big).
\end{align}
Here and below $C$ is a ``large'' constant depending on the value
of some ``known'' norms of the two solutions. Here, for instance, $C$
may depend on the $L^\infty(0,T;H)$-norm of $\sigma_2$ and on the 
$L^\infty(0,T;L^3(\Omega))$-norm of $\sigma_1$; the latter is controlled
due to \eqref{rego:sigma} and the continuous embedding 
$D(A^{1/4})\subset L^3(\Omega)$. 
%
 %
%

%
The estimation of the term depending on $k_1$ and $k_2$ is analogous and 
in fact simpler. Mimicking the above procedure we actually obtain
\begin{equation}\label{co:55}
  \bigg| \io \big( k_1 \sigma_1 - k_2 \sigma_2 \big) \bigg|
   \le c \big( 1 + \| \sigma_1 \| \big) \big( \| \sigma \| + \| \fhi \| \big)
   \le C \big( \| \sigma \| + \| \fhi \| \big).
\end{equation}
%
%
%
Hence, multiplying \eqref{co:53} by $2 \sigma\OO$, and  
using \eqref{co:54}-\eqref{co:55},
we infer
\begin{align}\no
  \frac12 \ddt |\sigma\OO|^2 
   & \le C |\sigma\OO| \big( 1 + \| \sigma_1 \|_{V} \big)
     \big( \| \sigma \| + \| \fhi \| \big)\\
 \label{diff:12}
   & \le C_\epsilon \big(1 + \| \sigma_1 \|_{V}^2 \big) | \sigma\OO |^2
     + \epsilon \| \fhi \|^2 + \epsilon \| \sigma - \sigma\OO \|^2,
\end{align} 
for ``small'' $\epsilon>0$ to be chosen later on and correspondingly 
``large'' $C_\epsilon>0$.

Next, we subtract \eqref{co:53} from \eqref{nutr:d} so to obtain
\begin{align}\no
  & (\sigma - \sigma\OO)_t 
   - \Delta\sigma + \chi \dive( \sigma_1 \nabla \fhi ) 
   + \chi \dive (\sigma \nabla \fhi_2)\\
 \no
  & \mbox{}~~~~~
   = \big( - h_1 \sigma_1^2 + h_2 \sigma_2^2 + k_1 \sigma_1 - k_2 \sigma_2 \big)
    - \big( - h_1 \sigma_1^2 + h_2 \sigma_2^2 + k_1 \sigma_1 - k_2 \sigma_2 \big)\OO\\
 \label{nutr:d2}
  & \mbox{}~~~~~
   =: H - H\OO.
\end{align} 
Testing \eqref{nutr:d2} by $\calN(\sigma - \sigma\OO)$ we deduce 
\begin{align}\no
  & \frac12 \ddt \| \sigma - \sigma\OO \|^2_*
   + \| \sigma - \sigma\OO \|^2 
   = \chi \io \sigma_1 \nabla\fhi \cdot \nabla \calN (\sigma-\sigma\OO)\\
 \label{co:56}
  & \mbox{}~~~~~
   + \chi \io \sigma \nabla\fhi_2 \cdot \nabla \calN (\sigma-\sigma\OO)
   + \io ( H - H\OO) \calN(\sigma-\sigma\OO).
\end{align}
In order to control the terms on the \rhs, we first observe that 
\begin{align}\no
  \chi \io \sigma_1 \nabla\fhi \cdot \nabla \calN (\sigma-\sigma\OO)
   & \le c \| \sigma_1 \|_{L^3(\Omega)}
    \| \nabla\fhi \| \|  \nabla \calN (\sigma-\sigma\OO) \|_{L^6(\Omega)}\\
  \label{co:57}
  & \le \frac16 \| \sigma - \sigma\OO \|^2
   + {\boldsymbol C} \| \nabla \fhi \|^2,
\end{align}
where the constant ${\boldsymbol C}$ depends
on the $L^\infty(0,T;L^3(\Omega))$-norm
of $\sigma_1$, which is a known quantity 
in view of \eqref{rego:sigma} and of the continuous 
embedding $D(A^{1/4})\subset L^3(\Omega)$.

Similarly, using also \eqref{nor:*2}, we have
\begin{align}\no
  \chi \io \sigma \nabla\fhi_2 \cdot \nabla \calN (\sigma-\sigma\OO)
   & \le c \| \sigma \|
    \| \nabla\fhi_2 \|_{L^\infty(\Omega)} \|  \nabla \calN (\sigma-\sigma\OO) \|\\
  \no
  & \le \frac1{12} \| \sigma \|^2
   + c \| \fhi_2 \|_{W^{2,6}(\Omega)}^2
    \| \sigma - \sigma\OO \|^2_*\\
  \label{co:58}
  & \le \frac16 \| \sigma - \sigma\OO \|^2
   + c | \sigma\OO |^2
   + c \| \fhi_2 \|_{W^{2,6}(\Omega)}^2
    \| \sigma - \sigma\OO \|^2_*.
\end{align}
It remains to control the last term on the \rhs\ of \eqref{co:56}. To 
this aim, we first point out that 
\begin{align}\no
  & \io (H-H\OO) \calN(\sigma-\sigma\OO)
    \le \| H-H\OO \|_{L^{6/5}(\Omega)} \| \calN(\sigma-\sigma\OO) \|_{L^6(\Omega)}\\
 \label{co:59}
  & \mbox{}~~~~~
   \le c \| H \|_{L^{6/5}(\Omega)} \| \calN(\sigma-\sigma\OO) \|_{V}
   \le c \| H \|_{L^{6/5}(\Omega)} \| \sigma-\sigma\OO \|_{*}.
\end{align}
To proceed, we observe that, by interpolation, one has
\begin{equation}\label{co:5a}
   \| \sigma_1 \|_V 
    \le c \| \sigma_1 \|_{D(A^{1/4})}^{1/2} \| \sigma_1 \|_{D(A^{3/4})}^{1/2}.
\end{equation}
Hence, using \eqref{rego:sigma}, it is easy to get
\begin{equation}\label{co:5b}
  \| \sigma_1 \|_{L^4(0,T;V)} \le C.
\end{equation}
To control the first factor on the \rhs\ of \eqref{co:59},
we mimick the procedure in~\eqref{co:54}-\eqref{co:55}. 
Namely, \eqref{co:54} is modified as follows:
\begin{align}\no
  & \| - h_1 \sigma_1^2 + h_2 \sigma_2^2 \|_{L^{6/5}(\Omega)}
   \le \| (h_1 - h_2) \sigma_1^2 \|_{L^{6/5}(\Omega)}
   + \| h_2 (\sigma_1 + \sigma_2) \sigma\|_{L^{6/5}(\Omega)}\\
 \no
   & \mbox{}~~~~~
    \le c \| h_1 - h_2 \| \| \sigma_1 \|_{L^6(\Omega)}^2
     + c \big( \| \sigma_1 \|_{L^3(\Omega)} + \| \sigma_2 \|_{L^3(\Omega)} \big) \| \sigma \|\\
 \no
   & \mbox{}~~~~~
    \le c \big( \| \sigma \| + \| \fhi \| \big) \| \sigma_1 \|_{L^6(\Omega)}^2
     + c \big( \| \sigma_1 \|_{L^3(\Omega)} + \| \sigma_2 \|_{L^3(\Omega)} \big) \| \sigma \|\\
 \label{co:54b}
   & \mbox{}~~~~~
    \le C \big( 1 +  \| \sigma_1 \|_{V}^2 \big) 
     \big( \| \sigma \| + \| \fhi \| \big).
\end{align}
Next, similarly with \eqref{co:55}, we have
\begin{equation}\label{co:55b}
  \| k_1 \sigma_1 - k_2 \sigma_2 \|_{L^{6/5}(\Omega)}
   \le C \big( \| \sigma \| + \| \fhi \| \big).
\end{equation}
Actually, we point out that
the constants $C$ on the \rhs s of \eqref{co:54b} and \eqref{co:55b}
also depend on the $L^\infty(0,T;L^3(\Omega))$-norms
of $\sigma_1$ and $\sigma_2$, which, again, are known quantities. Hence,
\eqref{co:59} yields
\begin{align}\no
  & \io (H-H\OO) \calN(\sigma-\sigma\OO) 
    \le C \big( 1 +  \| \sigma_1 \|_{V}^2 \big) 
     \big( \| \sigma - \sigma\OO \| + |\sigma\OO| + \| \fhi \| \big)
    \| \sigma-\sigma\OO \|_*\\
 \label{co:59b}
  & \mbox{}~~~~~
   \le C \big( 1 +  \| \sigma_1 \|_{V}^4 \big) \| \sigma-\sigma\OO \|_*^2
    + \frac16 \| \sigma - \sigma\OO \|^2  + c |\sigma\OO|^2 + c\| \fhi \|^2.
\end{align}
Now, taking \eqref{co:57}-\eqref{co:59b} into account, we deduce
from \eqref{co:56} that
\begin{align}\no
  & \frac12 \ddt \| \sigma - \sigma\OO \|^2_*
   + \frac12 \| \sigma - \sigma\OO \|^2 
   \le {\boldsymbol C} \| \fhi \|_V^2
   + c |\sigma\OO|^2\\
 \label{co:56b}
  & \mbox{}~~~~~
    + C \big( 1 + \| \sigma_1 \|_{V}^4 
        + \| \fhi_2 \|_{W^{2,6}(\Omega)}^2 \big) \| \sigma-\sigma\OO \|_*^2,
\end{align}
possibly for a new value of the constant ${\boldsymbol C}$ 
(which, we recall, depends on the $L^\infty(0,T;L^3(\Omega))$-norm 
of $\sigma_1$ and on the fixed parameters of the problem).

Let us now multiply \eqref{co:52} by $2{\boldsymbol C}$
and add the result to \eqref{co:56b} so to get
\begin{align}\no
  & \frac12 \ddt \| \sigma - \sigma\OO \|^2_*
   + 2{\boldsymbol C}\ddt \| \fhi \|_*^2
   + \frac12 \| \sigma - \sigma\OO \|^2 
   + {\boldsymbol C} \| \fhi \|_V^2\\
 \label{co:56c}
  & \mbox{}~~~~~
   \le c {\boldsymbol C} \| \fhi\|_*^2 
    + c |\sigma\OO|^2
    + C \big( 1 + \| \sigma_1 \|_{V}^4
        + \| \fhi_2 \|_{W^{2,6}(\Omega)}^2 + {\boldsymbol C} \big) 
             \| \sigma-\sigma\OO \|_*^2.
\end{align}
Then, we go back to \eqref{diff:12} and
choose $\epsilon=\min\{1/4,{\boldsymbol C}/2\}$ therein;
we then sum it to \eqref{co:56c} so to obtain
\begin{align}\no
  & \frac12 \ddt \| \sigma - \sigma\OO \|^2_*
   + \frac12 \ddt |\sigma\OO|^2 
   + 2{\boldsymbol C}\ddt \| \fhi \|_*^2
   + \frac14 \| \sigma - \sigma\OO \|^2 
   + \frac{{\boldsymbol C}}2 \| \fhi \|_V^2\\
 \label{co:56d}
  & \mbox{}~~~~~
   \le c {\boldsymbol C} \| \fhi\|_*^2 
    + C \big(1 + \| \sigma_1 \|_{V}^2 \big) | \sigma\OO |^2 
    + C \big( 1 + \| \sigma_1 \|_{V}^4 
        + \| \fhi_2 \|_{W^{2,6}(\Omega)}^2 + {\boldsymbol C} \big) 
            \| \sigma-\sigma\OO \|_*^2.
\end{align}
Then, using \eqref{co:5b} (for $\sigma_1$) and the last of 
\eqref{rego:fhi} (for $\fhi_2$)
and applying Gr\"onwall's lemma, we readily get the assert. 
The proof of Theorem~\ref{teo:wellpo} is complete.


\section{Asymptotic compactness and construction of the attractor}
\label{sec:asy}

In this part we detail the proof of Theorem~\ref{teo:attra}.
This will be achieved by showing additional regularity estimates
holding uniformly for sufficiently large values of the time variable.

To this purpose, we start recalling that,
for any set $B$ of initial data bounded in $\Phi_m\times \Sigma$,
any solution emanating from $B$
takes values in the bounded absorbing set $\calB_0$ for any 
$t\ge T_1$, where $T_1$ only depends on the ``radius'' of $B$ 
in $\Phi_m\times \Sigma$. With a simple argument based on the 
uniform Gr\"onwall lemma we will now see that any such 
solution takes values in a bounded set $\calB_1$ of 
$\Psi_m\times V_+$ for every $t\ge T_1 + 2$. Using the fact
that the immersion of $\Psi_m \times V_+$ 
into $\Phi_m \times \Sigma$ is compact (which is 
clear for what concerns $\sigma$ and follows from
the result \cite[Prop.~3.15]{RSS} reported here as 
Prop.~\ref{pro:RS} for what concerns $\fhi$),
the general theory of infinite-dimensional dynamical 
systems (applied to the {\sl closed}\/ semigroup $S(t)$,
cf.\ \cite{PZ}) then guarantees
the existence of the global attractor in the sense
of Theorem~\ref{teo:attra}.

To deduce the additional regularity estimates, we first
test \eqref{nutr} by $-\Delta\sigma$ so
to obtain
\begin{align}\no
  & \frac12 \ddt \| \nabla \sigma \|^2
   + \| \Delta\sigma \|^2
   = \chi \io \nabla\sigma \cdot \nabla \fhi \Delta \sigma\\
 \label{co:81}
  & \mbox{}~~~~~
   + \chi \io \sigma \Delta \fhi \Delta \sigma
   + \io \big( h(\sigma,\fhi) \sigma^2 - k(\sigma,\fhi) \sigma \big) \Delta\sigma.
\end{align}
Then, the terms on the \rhs\ are controlled as follows:
\begin{align}\label{co:82}
  & \chi \io \nabla\sigma \cdot \nabla \fhi \Delta \sigma
   \le c \| \nabla \sigma \| \| \nabla \fhi \|_{L^\infty(\Omega)} \| \Delta \sigma \|
   \le \frac14 \| \Delta \sigma \|^2
    + c \| \nabla \sigma \|^2 \| \fhi \|_{W^{2,6}(\Omega)}^2,\\
 \label{co:83}
  & \chi \io \sigma \Delta \fhi \Delta \sigma
   \le c \| \sigma \|_{L^6(\Omega)} \| \Delta \fhi \|_{L^3(\Omega)} \| \Delta \sigma \|
   \le \frac14 \| \Delta \sigma \|^2
    + c \| \sigma \|_V^2 \| \fhi \|_{W^{2,3}(\Omega)}^2,\\
 \label{co:84}
  & \io \big( h(\fhi,\sigma) \sigma^2 - k(\fhi,\sigma) \sigma \big) \Delta\sigma
   \le \frac14 \| \Delta \sigma \|^2
    + c \big( 1 + \| \sigma \|_{L^4(\Omega)}^4\big)
   \le \frac14 \| \Delta \sigma \|^2
    + c \big( 1 + \| \sigma \|_{V}^4\big).
\end{align}
Collecting the above computations, \eqref{co:81} gives
\begin{equation}\label{co:85}
  \frac12 \ddt \| \nabla \sigma \|^2
   + \frac14 \| \Delta\sigma \|^2
   \le c \big( 1 + \| \sigma \|_{V}^2 + \| \fhi \|_{W^{2,6}(\Omega)}^2 \big)
    \| \sigma \|_V^2.
\end{equation}
Hence, recalling \eqref{mo:27} and \eqref{mo:31}, and using the uniform Gr\"onwall 
lemma, we deduce 
\begin{equation}\label{co:86}
  \| \sigma(t) \|_{V}
   + \| \sigma \|_{L^2(t,t+1;W)}
   \le Q(R_0),  
\end{equation}
say, for every $t\ge T_1 + 1$.
At this point, a comparison of terms in \eqref{nutr} (equivalently, one could
also use $\sigma_t$ as a test function therein) also gives
\begin{equation}\label{co:87}
  \| \sigma_t \|_{L^2(t,t+1;H)} \le Q(R_0),  
\end{equation}
still for $t\ge T_1 + 1$.

With \eqref{co:86}-\eqref{co:87} at disposal, we can go back to 
the Cahn-Hilliard system \eqref{CH1}-\eqref{CH2} and perform the so-called
second-energy estimate. Namely, we test \eqref{CH1} by $\mu_t$ and 
the time derivative of \eqref{CH2} by $\fhi_t$ (also this argument can 
be easily justified in an approximation scheme).
Using the monotonicity of $\beta$, it is then easy to arrive at 
\begin{equation}\label{co:88}
  \frac12 \ddt \| \nabla \mu \|^2
   + \| \nabla\fhi_t \|^2
   \le \lambda \| \fhi_t \|^2 
   + \chi \io \sigma_t \fhi_t.
\end{equation}
Using Ehrling's lemma, the terms on the \rhs\ are controlled
as follows:
\begin{equation}\label{co:89}
  \lambda \| \fhi_t \|^2 
   + \chi \io \sigma_t \fhi_t
   \le \frac12 \| \nabla\fhi_t \|^2
   + c \| \fhi_t \|_*^2
   + c \| \sigma_t \|^2,
\end{equation}
where we also used the fact that $\fhi_t$ has zero spatial average due to
conservation of mass. Replacing \eqref{co:89} into \eqref{co:88}, we obtain
\begin{equation}\label{co:88b}
   \ddt \| \nabla \mu \|^2
   + \| \nabla\fhi_t \|^2
   \le  c \| \fhi_t \|_*^2
   + c \| \sigma_t \|^2.
\end{equation}
Hence, recalling \eqref{di:14}-\eqref{di:15} and \eqref{co:87}, we may apply 
once more the uniform Gr\"onwall lemma so to deduce 
\begin{equation}\label{co:90}
  \| \nabla \mu(t) \|^2
   + \int_t^{t+1} \| \nabla\fhi_t \|^2
   \le Q(R_0),
\end{equation}
for every $t\ge T_1+2$. More precisely, using \eqref{mz:11}
and \eqref{di:15}, the above relation is improved to
\begin{equation}\label{co:90b}
  \| \mu(t) \|^2_V
   + \int_t^{t+1} \| \fhi_t \|_V^2
   \le Q(R_0),
\end{equation}
still for $t\ge T_1+2$. Now, going back to inequality \eqref{mo:26}
and observing that its \rhs\ is now controlled uniformly in time thanks to 
\eqref{co:86} and \eqref{co:90b}, we readily deduce
\begin{equation}\label{co:91}
  \| \beta(\fhi(t)) \|_{L^6(\Omega)}
   + \| \fhi(t) \|_{W^{2,6}(\Omega)}
   \le Q(R_0),
\end{equation}
for $t\ge T_1+2$. Then, we observe that, for what concerns 
$\sigma$, the regularity prescribed in \eqref{bound:A}
corresponds exactly to the outcome of estimate
\eqref{co:86}, whereas, regarding $\fhi$, \eqref{co:91}
yields directly \eqref{bound:A1} (which is of course stronger
than \eqref{bound:A}). The proof of Theorem~\ref{teo:attra}
is then complete.
\beos\label{oss:beta0}
 As already observed in Remark~\ref{oss:multi}, the first condition 
 in \eqref{co:91} should be intended to hold for the section $\eta$
 of the multi-function $\beta(\fhi)$ satisfying \eqref{CH2-multi}.
 As the definitions \eqref{defi:Psim} of the space $\Psi_m$ and 
 \eqref{dist:Psim} of the related metric use the ``minimal section''
 $\beta^0$, it may be worth observing that, as $\beta^0$ is, 
 indeed, minimal, whenever the first \eqref{co:91} holds
 for some section $\eta$, it consequently holds also for $\beta^0(\fhi)$. 
\eddos

%
 %
 %
%
%
%


\subsection{Further regularity of the attractor}
\label{subsec:H2}

We prove here the additional asymptotic regularity of $\sigma$
stated in Theorem~\ref{teo:regattra}.
To this aim, rewriting equation \eqref{nutr} in the form 
\begin{equation}\label{nutr:A2}
  \sigma_t + A \sigma = - \chi \nabla \sigma\cdot \nabla\fhi
   - \chi \sigma \Delta \fhi
   - h(\sigma,\fhi) \sigma^2 + ( k(\sigma,\fhi) + 1) \sigma =: R(\sigma,\fhi),
\end{equation}
using properties \eqref{co:86} and \eqref{co:91} with Sobolev's embeddings,
it is not difficult to infer
\begin{equation}\label{H2:11}
  \| R(\sigma,\fhi) \|_{L^\infty(t,t+1;H)} \le Q(R_0),
\end{equation}
at least for $t\ge T_1 + 2$.
In addition to that, we may observe that, as a further consequence
of \eqref{co:86}, for every $t\ge T_1 + 2$ there exists ${\ov t}\in [t,t+1]$
such that 
\begin{equation}\label{H2:11b}
  \| \sigma({\ov t}) \|_{H^2(\Omega)} \le Q(R_0).
\end{equation}
Hence, applying over the time interval $[{\ov t},{\ov t}+2]$ 
parabolic regularity estimates of Agmon-Douglis-Nirenberg type
to problem \eqref{nutr:A2} with the ``initial'' condition 
$\sigma({\ov t})$, we readily deduce
\begin{equation}\label{H2:12x}
  \| \sigma_t \|_{L^p({\ov t},{\ov t}+2;H)}
   + \| A \sigma \|_{L^p({\ov t},{\ov t}+2;H)} 
   \le Q(R_0),
\end{equation}
which of course implies
\begin{equation}\label{H2:12}
  \| \sigma_t \|_{L^p(t,t+1;H)}
   + \| A \sigma \|_{L^p(t,t+1;H)}
   \le Q(R_0),
\end{equation}
at least for $t\ge T_1+3$. Note that \eqref{H2:12x} and \eqref{H2:12} 
hold for $p\in[1,\infty)$ with the expression of $Q$ depending on $p$
and possibly ``exploding'' as $p\nearrow\infty$.

Then, choosing $p$ large enough (but finite), applying standard
elliptic regularity and interpolation results, we also deduce
\begin{equation}\label{H2:13}
  \| \sigma \|_{L^\infty(t,t+1;L^\infty(\Omega))}
   \le Q(R_0).
\end{equation}
The above uniform boundedness property is the key step for 
obtaining additional estimates. Actually, we can now take the time derivative 
of \eqref{nutr} and test it by $\sigma_t$ so to deduce
\begin{equation}\label{H2:14}
  \frac12\ddt \| \sigma_t \|^2
   + \| \nabla\sigma_t \|^2
   = \chi \io \sigma_t \nabla\fhi \cdot \nabla\sigma_t
   + \chi \io \sigma \nabla\fhi_t \cdot \nabla\sigma_t
   + \io \big( - h(\sigma,\fhi) \sigma^2 
        + k(\sigma,\fhi) \sigma \big)_t \sigma_t.
\end{equation}
The terms on the \rhs\ can be controlled as follows:
\begin{align}\label{H2:15}
  & \chi \io \sigma_t \nabla\fhi \cdot \nabla\sigma_t
   \le \| \sigma_t \| \| \nabla \fhi \|_{L^\infty(\Omega)} \| \nabla\sigma_t \|
    \le \frac14 \| \nabla\sigma_t \|^2 + Q(R_0) \| \sigma_t \|^2,\\
 \label{H2:16}
  & \chi \io \sigma \nabla\fhi_t \cdot \nabla\sigma_t
   \le \| \sigma \|_{L^\infty(\Omega)}  \| \nabla \fhi_t \| \| \nabla\sigma_t \|
    \le \frac14 \| \nabla\sigma_t \|^2 + Q(R_0) \| \nabla\fhi_t \|^2,\\
 \label{H2:17}
  & \io \big( - h(\sigma,\fhi) \sigma^2 + k(\sigma,\fhi) \sigma \big)_t \sigma_t
   \le Q(R_0) \big( \| \sigma_t \|^2 + \| \fhi_t \|^2 \big),
\end{align}
where we have also used \eqref{co:91} with the continuous embedding 
$W^{2,6}(\Omega)\subset W^{1,\infty}(\Omega)$, \eqref{H2:13},
and Assumption~(A1).
 
Using \eqref{H2:15}-\eqref{H2:17}, \eqref{H2:14} gives
\begin{equation}\label{H2:14b}
  \ddt \| \sigma_t \|^2
   + \| \nabla\sigma_t \|^2
   \le Q(R_0) \big( \| \sigma_t \|^2 + \| \fhi_t \|_V^2 \big),
\end{equation}
whence, using \eqref{co:87}, \eqref{co:90b} and the uniform
Gr\"onwall lemma, we deduce 
\begin{equation}\label{H2:18}
  \| \sigma_t(t) \|
   + \| \sigma_t \|_{L^2(t,t+1;V)} \le Q(R_0),
\end{equation}
for every $t\ge T_1 + 4$. Finally, interpreting \eqref{nutr} 
as a time-dependent family of elliptic problems, namely
\begin{equation}\label{nutr:ell}
   A \sigma = - \sigma_t 
   - \chi \nabla\sigma\cdot\nabla\fhi 
   - \chi \sigma \Delta\fhi 
   - h(\sigma,\fhi) \sigma^2 + ( k(\sigma,\fhi) + 1) \sigma,
\end{equation}
we may observe that, as a consequence of \eqref{co:91} and \eqref{H2:18},
the \rhs\ is bounded in $H$ uniformly for $t\ge T_1 + 4$.
Hence, applying once more standard elliptic regularity results, we deduce
\begin{equation}\label{H2:19}
  \| \sigma(t) \|_{H^2(\Omega)} \le Q(R_0),
\end{equation}
for every $t\ge T_1 + 4$. As this corresponds exactly to 
\eqref{bound:A2}, the proof of Theorem~\ref{teo:regattra} is complete.


\section{Proof of Theorem~\ref{teo:ln}}
\label{subsec:sig}

In this section we analyze in more detail the sign properties of 
$\sigma$, so to prove Theorem~\ref{teo:ln}. To this aim, partially
following~\cite{LRS}, we first provide a sort of ``entropic''
reformulation of \eqref{nutr}, which is obtained by 
multiplying \eqref{nutr} by $\sigma^{-1}$ and performing some 
integrations by parts. We point out that also this argument
is formal; indeed, at this level we do not know
whether $\sigma^{-1}$ has any summability property and
in principle we could not even exclude that there may exist
regions of positive measure where $\sigma$ is 
identically $0$ (we will know that this in fact cannot happen
as {\sl an outcome}\/ of the procedure). A simple way to make 
our procedure rigorous could be that of multiplying \eqref{nutr}
by $T_n(\sigma^{-1})$, where $T_n$ is the truncation operator 
at height $n$, and working on the relation obtained that way. 
Then, as the resulting estimates will turn out to be 
independent of $n$, at the end one may let $n\to \infty$. 
For the sake of brevity we omit the details of this argument
and work directly on relation \eqref{nutr}, whose ``entropic''
reformulation then reads
\begin{equation}\label{nutr:v}
  v_t - \Delta v - | \nabla v |^2
   + h(\sigma,\fhi) \sigma - k(\sigma,\fhi)
   = - \chi \nabla v \cdot \nabla \fhi
   - \chi \Delta \fhi,
\end{equation} 
and where we have set $v:=\ln \sigma$. Let us now consider
a generic monotone function $\gamma:\RR\to\RR$, let 
$\gammaciapo$ be some primitive of $\gamma$, and 
let us test \eqref{nutr:v} by $\gamma(v)$ so to deduce 
the general formula
\begin{align}\no
  &  \ddt \io \gammaciapo(v) 
   + \io \big( \gamma'(v) - \gamma(v) \big) | \nabla v |^2
   - \io k(\sigma,\fhi) \gamma(v) 
   + \io h(\sigma,\fhi) \sigma \gamma(v) \\
 \label{vv:11}  
  & \mbox{}~~~~~
  = - \chi \io \gamma(v) \nabla v \cdot \nabla \fhi
   - \chi \io \gamma(v) \Delta \fhi.
\end{align} 
Now, as we need to estimate the values of $\sigma$ that are close 
to $0$, we start with taking $\gamma(v) = -1$ for $v< 0$ and 
$\gamma(v)=0$ for $v\ge 0$; namely, 
$\gamma(v) = -(\textrm{sign}(v))_{-}= \min\left\{\textrm{sign}(v), 0\right\}$ in 
such a way that 
$\gammaciapo(v) = v_{-}$.
Of course $\gamma$ is a monotone
function, but it is discontinuous at $0$, hence the procedure 
is formal; to make it rigorous one could replace $\gamma$ 
by a function $\gamma_\delta$ taking the value $\delta^{-1}r$
in the interval $[-\delta,0]$ and coinciding with $\gamma$
elsewhere, and then let $\delta\searrow 0$. We omit the details
since also this argument is standard. 
Then, recalling that, if $v\in V$, then $\nabla v_{-}\in H$ with 
\begin{equation}\label{as:11}  
\nabla v_{-}
=
\begin{cases}
-\nabla v\qquad &\textrm{ if }\, v\le 0\\
0\qquad &\textrm{ otherwise }
\end{cases}
\end{equation}
almost everywhere in $\Omega$ (thanks to \cite[Lemme 1.1 \& Lemme 1.2]{Stampa}
or \cite[Lemma 7.6]{GT})
and neglecting some nonnegative terms from the \lhs, \eqref{vv:11}
reduces to 
\begin{equation}\label{vv:12}  
  \ddt \| v_- \|_{L^1(\Omega)}
   + \| \nabla v_- \|^2
   \le \io h(\sigma,\fhi) \sigma 
   + \chi \io | \nabla v_-| | \nabla \fhi |
   + \chi \io | \Delta \fhi |.
\end{equation} 
Now, integrating the above over $(0,T)$ and using \eqref{rego:fhi}
and \eqref{rego:sigma}, it is a standard matter to control the terms on the 
\rhs. Hence, using also condition \eqref{strongpos} on the initial datum, we
deduce 
\begin{equation}\label{stv:1}  
  \| v_- \|_{L^\infty(0,T;L^1(\Omega))}
   + \| v_- \|_{L^2(0,T;V)}
  \le Q(T),
\end{equation} 
where the \rhs\ depends on the reference
time $T$ as the above estimate does not have a dissipative 
character (we refer to Remark~\ref{oss:nodiss} for additional 
comments).

Notice also that \eqref{stv:1}, combined with Sobolev's embeddings,
implies in particular that, for any ``small'' $\tau\in(0,T)$, 
there exists some $t_\tau \in (\tau/2,\tau)$  such that 
\begin{equation}\label{stv:2}
  \| v_-(t_\tau) \|_{L^6(\Omega)} \le c \| v_-(t_\tau) \|_{V} 
   \le Q(T,\tau^{-1}).
\end{equation} 
Interpreting the above quantity as an ``initial'' datum for 
equation \eqref{nutr:v}, we will now derive
additional estimates holding for strictly positive times. To this aim, we
repeat \eqref{vv:11}, taking now 
$\gamma(v)=-(v_-)^{p-1}$ for a generic exponent $p\ge 2$. Then, we readily 
obtain
\begin{align}\no
  & \frac1{p} \ddt \| v_- \|_{L^p(\Omega)}^p
   + \io \big( v_-^{p-1} + (p-1) v_-^{p-2} \big) | \nabla v_- |^2
   + \io k(\sigma,\fhi) v_-^{p-1}
   - \io h(\sigma,\fhi) \sigma v_-^{p-1}\\
 \label{vv:13}  
  & \mbox{}~~~~~
   \le \chi \io v_-^{p-1} \nabla v \cdot \nabla \fhi
    + \chi \io v_-^{p-1} \Delta \fhi
   = \chi \io v_-^{p-1} \nabla v_- \cdot \nabla \fhi
    + \chi \io v_-^{p-1} \Delta \fhi.
\end{align} 
In order to deal with the above relation, we go back to 
estimate~\eqref{co:91} and make the following general consideration:
modifying a bit the procedure given in the previous section and, in particular, 
applying Gr\"onwall lemma on intervals of length $\tau$ (or, more precisely,
$\tau/n$ for sufficiently large $n\in\NN$, where $n$ is the number
of iterations), rather than on intervals of length $1$,
one could restate \eqref{co:91} as 
\begin{equation}\label{co:91tau}
  \| \beta(\fhi(t)) \|_{L^6(\Omega)}
   + \| \fhi(t) \|_{W^{2,6}(\Omega)}
   \le Q(\tau^{-1}),
   \quext{for every }\,t\in[\tau,T]
    ~~\text{and every }\,\tau>0,
\end{equation}
where now the expression of $Q$ may also depend on the specific choice 
of the initial datum as we are not looking for a dissipative estimate 
at this level.
Hence, using \eqref{co:91tau} with Young's inequality,
there follows
\begin{equation}\label{vv:14}
  \chi \io v_-^{p-1} \nabla v_- \cdot \nabla \fhi
   \le \frac12 \io v_-^{p-1} | \nabla v_- |^2
    + Q(\tau^{-1}) \| v_- \|_{L^{p-1}(\Omega)}^{p-1}.
\end{equation} 
Let us now consider the last two terms on the \lhs. To this aim, 
let us observe that, applying the above considerations to 
\eqref{H2:19} and using Sobolev's embeddings, there follows
\begin{equation}\label{H2:19b}
  \| \sigma \|_{L^\infty(\tau,T;L^\infty(\Omega))}
  \le c \| \sigma \|_{L^\infty(\tau,T;H^2(\Omega))}
   \le Q(\tau^{-1}).
\end{equation}
Using \eqref{H2:19b} and Assumption~(A1), it is then clear that 
\begin{equation}\label{vv:15}
  \bigg| \io \big( k(\sigma,\fhi) - h(\sigma,\fhi) \sigma \big) v_-^{p-1} \bigg|
   \le Q(\tau^{-1}) \| v_- \|_{L^{p-1}(\Omega)}^{p-1}. 
\end{equation}
Next, we observe that the second term on the \lhs\ of \eqref{vv:13}
(where we reduced the value of one constant as we have 
incorporated there the contribution of the first term on 
the \rhs\ of \eqref{vv:14}) can be equivalently rewritten as 
\begin{equation}\label{vv:16}
  \io \Big( \frac12 v_-^{p-1} + (p-1) v_-^{p-2} \Big) | \nabla v_- |^2
   = \frac{2}{(p+1)^2} \big\| \nabla v_-^{\frac{p+1}2} \big\|^2
   + \frac{4(p-1)}{p^2} \big\| \nabla v_-^{\frac{p}2} \big\|^2.
\end{equation}
Then, we control the last term in \eqref{vv:13} as follows:
\begin{align}\no
  \chi \io v_-^{p-1} \Delta \fhi
   & \le c \| v_-^{\frac{p-2}2} \|
     \| v_-^{\frac{p}2} \|_V \| \Delta \fhi \|_{L^3(\Omega)}\\
 \no
   & \le Q(\tau^{-1}) \| v_-^{\frac{p-2}2} \|
      \| v_-^{\frac{p}2} \|_V \\
 \no
   & \le Q(\tau^{-1}) \| v_-^{\frac{p-2}2} \|
     \big( \| \nabla v_-^{\frac{p}2} \| + \| v_-^{\frac{p}2} \|_{L^1(\Omega)} \big) \\
 \label{vv:17}    
   & \le \frac{2(p-1)}{p^2} \| \nabla v_-^{\frac{p}2} \|^2
    + Q(\tau^{-1}) \big( 1 + p \| v_- \|_{L^{p-1}(\Omega)}^{p-1} \big),
\end{align}
where we have also used Young's inequality and the fact $p/2 \le p-1$ 
as we have assumed $p\ge 2$. Collecting 
\eqref{vv:14}-\eqref{vv:17} then \eqref{vv:13} gives
\begin{align}\no
  & \frac1{p} \ddt \| v_- \|_{L^p(\Omega)}^p
   + \frac{2}{(p+1)^2} \big\| \nabla v_-^{\frac{p+1}2} \big\|^2
   + \frac{2(p-1)}{p^2} \big\| \nabla v_-^{\frac{p}2} \big\|^2\\
 \label{vv:13b}
  & \mbox{}~~~~~
   \le Q(\tau^{-1}) \big( 1 + p \| v_- \|_{L^{p-1}(\Omega)}^{p-1} \big).
\end{align} 
The above differential inequality may serve as a starting point for 
a Moser iteration scheme with regularization, and, more precisely,
we may follow the lines of the procedure given in
\cite[Proof of Lemma~3.3]{SSZ} (see in particular the differential 
inequality~\cite[(3.52)]{SSZ}, which has exactly the same structure
of or \eqref{vv:13b}), or the similar argument in 
\cite[End of Proof of Thm 2.2]{iopf} (actually the situation
here is much simpler, in the sense that the regularity of the
forcing terms in \eqref{vv:12} is higher compared
to both the quoted references). 

We omit to report the details of the Moser scheme, as 
they are rather technical: the spirit of the procedure is to
perform infinitely many iterations by working on time intervals of
smaller and smaller length behaving like $\tau 2^{-k}$, $k\in \NN$,
and exploiting at each step the the parabolic
regularization effects. At the first iteration we can 
use \eqref{stv:2} as a regularity condition on the ``initial''
datum; then \eqref{stv:2} will be improved step by step
exploiting the effects of the diffusion term.

It can then be shown that the resulting bounds are uniform 
for large $p$ and that the behavior with respect to $\tau$ can also be 
quantitatively controlled (see the statement of 
\cite[Lemma~3.3]{SSZ}, where the argument was carried out for $\tau=1$, 
and in particular compare formulas~(3.50) and (3.51) therein).
As a final outcome of the procedure we then get
\begin{equation}\label{vv:21}
  \| v_-(t) \|_{L^\infty(\Omega)}
   \le Q(\tau^{-1},T), \quext{for every }\,t\ge \tau>0,
\end{equation}
which clearly implies \eqref{strongpos3} as the logarithm is 
unbounded near $0$.
\beos\label{vers:u}
 We observe that one could obtain a variant of the above result by
 an alternative procedure. Namely, if we multiply \eqref{nutr} 
 by $-1/\sigma^2$ we deduce
 \begin{equation}\label{nutr:2}
   \ddt \frac{1}{\sigma}
    + \frac 1{\sigma^2}\dive(\nabla\sigma)
    -  \chi \frac {\nabla\sigma}{\sigma^2} \cdot \nabla \fhi
    -  \chi \frac1{\sigma} \Delta \fhi 
    = h(\sigma,\fhi) - \frac{k(\sigma,\fhi)}{\sigma}.
 \end{equation}  
 Then, noting that 
 \begin{equation}\label{nutr:3}
   \frac 1{\sigma^2}\dive(\nabla\sigma)
    = \dive\Big(\frac{\nabla\sigma}{\sigma^2}\Big)
     - \nabla\sigma \cdot \nabla\Big(\frac1{\sigma^2}\Big)
    = - \Delta\Big(\frac1\sigma\Big) + 2 \frac{|\nabla\sigma|^2}{\sigma^3}
    = - \Delta\Big(\frac1\sigma\Big) + 8 |\nabla\sigma^{-1/2}|^2,
 \end{equation} 
 and setting $u:=\sigma^{-1}$, we arrive at the relation 
 \begin{equation}\label{nutr:u}
   u_t - \Delta u + 8 | \nabla u^{1/2} |^2
    + \chi \nabla u \cdot \nabla \fhi 
    - \chi u \Delta \fhi
   = h(\sigma,\fhi) - k(\sigma,\fhi) u,
 \end{equation} 
 which may also serve as a starting point for an alternative
 version of the Moser iteration
 argument. We preferred, however, to start from \eqref{nutr:v}
 as this requires the weaker assumption \eqref{strongpos}
 on the initial datum, whereas using \eqref{nutr:u} we should likely
 assume some summability on $u_0=\sigma_0^{-1}$, which is of 
 course a more restrictive condition.
\eddos


\section*{Acknowledgments}

G.S.~has been
 partially supported by the PRIN MUR Grant 2020F3NCPX 
``Mathematics for industry 4.0 (Math4I4)'' and by GNAMPA (Gruppo Nazionale per l'Analisi 
Matematica, la Probabilit\`a e le loro Applicazioni) of INdAM (Istituto Nazionale
di Alta Matematica). 
A.S.~has been partially supported by PRIN 2022 (Project no.\ 2022J4FYNJ), funded by MUR, Italy,
and the European Union -- Next Generation EU, Mission 4 Component 1 CUP F53D23002760006 and by 
GNAMPA (Gruppo Nazionale per l'Analisi 
Matematica, la Probabilit\`a e le loro Applicazioni) of INdAM (Istituto Nazionale
di Alta Matematica).



\vspace{15mm}

\noindent%
{\bf First author's address:}\\[1mm]
Giulio Schimperna\\
Dipartimento di Matematica, Universit\`a degli Studi di Pavia\\
and Istituto di Matematica Applicata e Tecnologie Informatiche ``Enrico Magenes'' (IMATI),\\
Via Ferrata, 5,~~I-27100 Pavia,~~Italy\\
E-mail:~~{\tt giulio.schimperna@unipv.it}

\medskip
\bigskip

\noindent%
{\bf Second author's address:}\\[1mm]
Antonio Segatti\\
Dipartimento di Matematica, Universit\`a degli Studi di Pavia\\
and Istituto di Matematica Applicata e Tecnologie Informatiche ``Enrico Magenes'' (IMATI),\\
Via Ferrata, 5,~~I-27100 Pavia,~~Italy\\
E-mail:~~{\tt antonio.segatti@unipv.it}


\begin{thebibliography}{99}




\bibitem{AS}
 A.~Agosti and A.~Signori, 
 {\sl Analysis of a multi-species Cahn-Hilliard-Keller-Segel tumor growth model with 
   chemotaxis and angiogenesis},
 J.~Differential Equations,
 {\bf 403} (2024), 
 308--367.

 
\bibitem{Ba}
 V.~Barbu,
 Nonlinear Semigroups and Differential Equations in Banach Spaces.
 Noordhoff, Leiden, 1976.

 
\bibitem{Br}
 H.~Br\'ezis,
 Op\'erateurs Maximaux Monotones et S\'emi-groupes de Contractions
    dans les Espaces de Hilbert.
 North-Holland Math.\ Studies {\bf 5},
 North-Holland,
 Amsterdam,
 1973.

 
\bibitem{CGSS}
 P.~Colli, G.~Gilardi, A.~Signori, and J.~Sprekels, 
 {\sl Solvability and optimal control of a multi-species 
   Cahn-Hilliard-Keller-Segel tumor growth model},
 ESAIM Control Optim.\ Calc.\ Var.,
 {\bf 31} (2025), Paper No. 85.
 
 
 \bibitem{CL}
 V.~Cristini and J.S.~Lowengrub,
 Multiscale Modeling of Cancer: An Integrated Experimental and Mathematical Modeling Approach. 
 Cambridge University Press, Leiden, 2010.
 
 
 \bibitem{CLLW}
 V.~Cristini, X.~Li, J.S.~Lowengrub, and S.M.~Wise,
 {\sl Nonlinear simulations of solid tumor growth using a mixture model: invasion 
  and branching},
 J.~Math.\ Biol.,
 {\bf 58} (2009), 723--763.
 
  
\bibitem{FLRS}
 S.~Frigeri, K.F.~Lam, E.~Rocca, and G.~Schimperna, 
 {\sl On a multi-species Cahn-Hilliard-Darcy tumor growth model with singular potentials},
 Commun.\ Math.\ Sci., {\bf 16} (2018), 
 821-856.
 
 
\bibitem{GL1}
 H.~Garcke and K.F.~Lam, 
 {\sl Global weak solutions and asymptotic limits of a Cahn-Hilliard-Darcy system 
   modelling tumour growth},
 AIMS Math.,
 {\bf 1} (2016), 318--360.
 
 
\bibitem{GL2}
 H.~Garcke and K.F.~Lam, 
 {\sl Well-posedness of a Cahn–Hilliard system modelling tumour growth with
   chemotaxis and active transport},
 Eur.\ J.~Appl.\ Math.,
 {\bf 28} (2017), 
 284--316.
 
 
\bibitem{GSS}
 G.~Gilardi, A.~Signori, and J.~Sprekels, 
 {\sl Nutrient control for a viscous Cahn-Hilliard-Keller-Segel model with logistic source describing tumor 
   growth},
 Discrete Contin.\ Dyn.\ Syst.\ Ser.~S, {\bf 16} (2023), 
 3552--3572.
 
\bibitem{GT}
D.~Gilbarg, and N.~Trudinger, 
Elliptic Partial Differential Equations of Second Order, 
Classics in Mathematics, Springer-Verlag, Berlin, (2001). 
 
 
\bibitem{GHW}
 A.~Giorgini, J.~He, and H.~Wu,
 {\sl Global weak solutions to a Navier-Stokes-Cahn-Hilliard System with chemotaxis
   and mass transport: cross diffusion versus logistic degradation},
 arXiv:2412.05751 (2024).

 
\bibitem{GP}
 M.~Grasselli and V.~Pata, 
 {\sl Existence of a universal attractor for a fully hyperbolic 
   phase-field system},
 J.~Evol.\ Equ., {\bf 4} (2004), 
 27--51.
 
 
\bibitem{Ha}
 J.K.~Hale,
 Asymptotic behavior of dissipative systems.
 Math.\ Surveys Monogr., 25.
 American Mathematical Society, Providence, RI, 
 1988.
 
 
\bibitem{LRS}
 R.~Lasarzik, E.~Rocca, and G.~Schimperna, 
 {\sl Weak solutions and weak-strong uniqueness for a Cahn-Hilliard type model with chemotaxis},
 paper in preparation. 

 
\bibitem{MiZe}
 A.~Miranville and S.~Zelik,
 {\sl Robust exponential attractors for Cahn-Hilliard type equations with singular potentials},
 Math.\ Methods Appl.\ Sci.,
 {\bf 27} (2004), 
 545--582.

 
 
 
\bibitem{PZ}
 V.~Pata and S.~Zelik,
 {\sl A result on the existence of global attractors for semigroups of closed operators},
 Commun.\ Pure Appl.\ Anal.,
 {\bf 6} (2007), 
 481--486.
 
 
\bibitem{RS}
 E.~Rocca and G.~Schimperna, 
 {\sl Universal attractor for some singular phase transition systems},
 Physica D: Nonlinear Phenomena,
{\bf 192} (2004), 279--307.

 
\bibitem{RSS}
 E.~Rocca, G.~Schimperna, and A.~Signori,
 {\sl On a Cahn-Hilliard-Keller-Segel model with generalized logistic source describing tumor growth},
 J.~Differential Equations,
 {\bf 343} (2023), 530--578.
 

\bibitem{iopf}
 G.~Schimperna,
 {\sl Global and exponential attractors for the Penrose-Fife system},
 Math.\ Models Methods Appl.\ Sci.,
 {\bf 19} (2009), 
 969--991.


\bibitem{Ssens}
 G.~Schimperna, 
 {\sl On a modified Cahn-Hilliard-Brinkman model with chemotaxis and nonlinear sensitivity},
 arXiv:2411.12505.
 
 
\bibitem{SSZ}
 G.~Schimperna, A.~Segatti, and S.~Zelik,
 {\sl Asymptotic uniform boundedness of energy solutions to the 
   Penrose-Fife model},
 J.~Evol.\ Equ., {\bf 12} (2012), 863--890.
 
 
 

\bibitem{Stampa}
G.~Stampacchia,
{\sl Le problème de Dirichlet pour les équations elliptiques du second ordre à coefficients discontinus},
Ann. Ist. Fourier (Grenoble), 
{\bf 15} (1965),
189--258
 
\bibitem{Te}
 R.~Temam, 
 Infinite-dimensional Dynamical Systems in Mechanics and Physics.
 Springer-Verlag, New York, 1988.
 
 
 
 
\bibitem{WLFC} 
 S.M.~Wise, J.S.~Lowengrub, H.B.~Frieboes, and V.~Cristini,
 {\sl Three-dimensional multispecies nonlinear tumor growth -- I: 
  model and numerical method},
 J.~Theor.\ Biol.,
 {\bf 253} (2008), 
 524--543.
 
 
  
\end{thebibliography}
\end{document}